\documentclass[11pt,leqno,draft]{amsart} 
\begin{document}

\author[C. Simpson]{Carlos Simpson}
\address{CNRS, Laboratoire J. A. Dieudonn\'e, UMR 6621
\\ Universit\'e de Nice-Sophia Antipolis\\
06108 Nice, Cedex 2, France}
\email{carlos@math.unice.fr}
\urladdr{http://math.unice.fr/$\sim$carlos/} 

\title[Gabriel-Zisman files]{Files for Gabriel-Zisman localization}

\maketitle

\tableofcontents

We present simplified versions of the main files for Gabriel-Zisman localization (the only file which isn't
included is {\tt updateA.v} which contains things which should have gone in earlier parts of the set theory and category theory developments). 
All proofs are
taken out, and also all comments have been removed. For proofs as well as the comments in their original form, see the actual Coq files
which come with this posting (go to ``Other formats'' then ``Source'' on the abstract page). For a general discussion
of the historical and little mathematical issues, for the bibliography, and for a more specific discussion of the present proof
files, see the companion  preprint ``Explaining Gabriel-Zisman localization to the computer''. The document below was generated by using
{\tt coqdoc} to create html files, cutting and pasting those together and editing the resulting text file, then using 
the ``verbatim'' environment for each module.

\section{Instructions}

We recall here the basic instructions. Everything should compile under the current version of Coq, version 8.0 patch level 2.  
See 
\verb}http://coq.inria.fr/distrib-eng.html} for downloading and installing
the Coq proof assistant.

The following instructions are for a Unix-like system. Installation, compilation and utilisation of large multiple-file projects in the 
Windows version of Coq has, for me at least, been problematic (it would be good if this could be improved in upcoming
versions of Coq). 

Go to the {\tt http://arxiv.org} website; 
then to the
abstract page for the present preprint. From the abstract page, click on ``Other formats'' 
for other download types
(rather than ``ps'' or ``pdf''). At the bottom of the resulting page under the rubrique
``Source'', click on ``Download source''. This yields a 
tar arxive called \verb}0506???.tar.gz} (fill in the correct paper number); 
save it in the directory you wish to use, then decompact it 
using the commands
\begin{verbatim}
gunzip 0506???.tar.gz
tar -xvf 0506???.tar
\end{verbatim}
The next steps are reproduced here from the instructions in the Coq reference manual.
First, create the makefile (if it isn't already present).
\begin{verbatim}
coq_makefile -o Makefile *.v
\end{verbatim}
Now create a file called \verb}.depend} and write into it the inter-file dependencies using the commands
\begin{verbatim}
touch .depend
make depend
\end{verbatim}
Subsequently if new files are added to the development, one should redo the command \verb}make depend}.

Now we are ready to compile everything! Type 
\begin{verbatim}
make all
\end{verbatim}
and you should get (quickly for the first two, then more slowly for the later files) the following list of compilations:
\begin{verbatim}
coqc  -q  -I .   axioms
coqc  -q  -I .   tactics
coqc  -q  -I .   set_theory
coqc  -q  -I .   functions
coqc  -q  -I .   notation
coqc  -q  -I .   universe
coqc  -q  -I .   order
coqc  -q  -I .   transfinite
coqc  -q  -I .   ordinal
coqc  -q  -I .   cardinal
coqc  -q  -I .   category
coqc  -q  -I .   functor
coqc  -q  -I .   nat_trans
coqc  -q  -I .   cat_examples
coqc  -q  -I .   limits
coqc  -q  -I .   colimits
coqc  -q  -I .   functor_cat
coqc  -q  -I .   fc_limits
coqc  -q  -I .   little_cat
coqc  -q  -I .   equalizer
coqc  -q  -I .   fiprod
coqc  -q  -I .   updateA
coqc  -q  -I .   freecat
coqc  -q  -I .   qcat
coqc  -q  -I .   gzdef
coqc  -q  -I .   left_fractions
coqc  -q  -I .   gzloc
coqc  -q  -I .   infinite
coqc  -q  -I .   lfcx
\end{verbatim}
The directory now contains both the \verb}v} files and the \verb}vo} (i.e. compiled) files. 
The files can be usefully perused using 
\verb}coqide}. Particularly useful is the \verb}Queries} set of commands in \verb}coqide}, allowing one to view 
important information from various files in useful chunks. 

The files which are new to the Gabriel-Zisman localization are the following:
\begin{verbatim}
updateA.v 
freecat.v 
qcat.v 
gzdef.v 
left_fractions.v  
gzloc.v 
infinite.v 
lfcx.v    
\end{verbatim}
The other files are identical to the set-theory and category-theory development that was bundled with 
the previous preprint math.AG/0410224.

The remainder of this paper is a listing of everything in the Gabriel-Zisman source files, except all the
proofs and the totality of the first file {\tt updateA.v}. This first file contains some things which logically
would go with the earlier set-theory and category-theory files; we group them into a new ``update'' file in order
to keep the earlier files identical to the previous version. 

For looking at the proofs it is better to use {\tt coqide}, the Coq integrated development
environment (which comes with Coq). For example to look at the file {\tt freecat.v}, do:
\begin{verbatim}
coqide freecat.v
\end{verbatim}
(or simply type {\tt coqide} then open the file you want). In {\tt coqide} you get a screen with the source
file
on the left, and the upper screen on the right is the proof goal window. Thus you can put the
cursor in the middle of a proof, click on the ``go here'' button, and see the state of the proof at that
point. Step through proofs with the ``next step'' button. Using these buttons requires having at 
least all of the previous files 
being compiled as per the instructions above---when {\tt coqide} gets to a {\tt Require Export}
instruction typically at the start of the file, it looks for a compiled {\tt vo} file.  

\newpage

\section{The file {\tt freecat.v}}

\begin{verbatim}

Library freecat
Set Implicit Arguments.
Unset Strict Implicit.

Require Export updateA.

\end{verbatim}

\subsection{The module {\tt Uple}}

\begin{verbatim}
Module Uple.
Export UpdateA.
Export Universe.

Definition axioms a :=
Function.axioms a &
inc (domain a) nat.

Definition length a := Bnat (domain a).

Definition create (f:nat -> E) (l:nat) :=
Function.create (R l) (fun x => f (Bnat x)).

Lemma domain_create : forall f l, domain (create f l) = R l.

Lemma length_create : forall f l, length (create f l) = l.

Definition component (i:nat) a :=
V (R i) a.

Lemma component_create : forall f l i,
i < l -> component i (create f l) = f i.

Lemma create_axioms : forall f l, axioms (create f l).

Lemma eq_create : forall a, axioms a ->
a = create (fun i => component i a) (length a).

Lemma compare_create : forall f g l m,
l = m -> (forall i, Peano.lt i l -> f i = g i) ->
create f l = create g m.

Lemma domain_emptyset : domain emptyset = emptyset.

Lemma empty_axioms : axioms emptyset.

Lemma length_emptyset : length emptyset = 0.

Lemma create_length_zero : forall f ,
create f 0 = emptyset.

Lemma uple_extensionality : forall a b,
axioms a -> axioms b ->
length a = length b ->
(forall i, i < (length a) -> component i a = component i b) ->
a=b.

Definition concatenate a b :=
create (fun i =>
Y (i < (length a)) (component i a)
(component (i - (length a)) b))
((length a) + (length b)).

Lemma length_concatenate : forall a b,
length (concatenate a b) = (length a) + (length b).

Lemma concatenate_axioms : forall a b,
axioms (concatenate a b).

Lemma component_concatenate_first : forall a b i,
i < (length a) -> component i (concatenate a b) = component i a.

Lemma component_concatenate_second : forall a b i,
i < (plus (length a) (length b)) ->
(length a) <= i ->
component i (concatenate a b) = component (i - (length a)) b.

Lemma component_concatenate_plus : forall a b i,
i < (length b) ->
component ((length a) + i ) (concatenate a b) = component i b.

Lemma concatenate_emptyset_left : forall a,
axioms a ->
concatenate emptyset a = a.

Lemma concatenate_emptyset_right : forall a,
axioms a ->
concatenate a emptyset = a.

Lemma concatenate_assoc : forall a b c,
concatenate a (concatenate b c) =
concatenate (concatenate a b) c.

Definition uple1 x :=
create (fun i:nat => x) 1.

Lemma uple1_axioms : forall x,
axioms (uple1 x).

Lemma length_uple1 : forall x, length (uple1 x) = 1.

Lemma component_uple1 : forall x i,
i < 1 -> component i (uple1 x) = x.

Lemma eq_uple1 : forall a x,
axioms a -> length a = 1 -> x = component 0 a ->
a = uple1 x.

Definition utack x a :=
concatenate a (uple1 x).

Lemma utack_axioms : forall a x,
axioms a -> axioms (utack x a).

Lemma length_utack : forall a x,
axioms a -> length (utack x a) = (length a) + 1.

Lemma domain_R_length : forall a,
axioms a -> domain a = R (length a).

Lemma inc_R_domain : forall a i,
axioms a ->
inc (R i) (domain a) = (i < (length a)).

Lemma component_utack_old : forall a x i,
axioms a -> i < (length a) ->
component i (utack x a) = component i a.

Lemma component_utack_new : forall a x i,
axioms a -> i = length a ->
component i (utack x a) = x.

Definition restrict a i :=
create (fun j => component j a) i.

Lemma length_restrict : forall a i,
length (restrict a i) = i.

Lemma component_restrict : forall a i j,
j < i -> component j (restrict a i) = component j a.

Lemma restrict_axioms : forall a i,
axioms (restrict a i).

Lemma eq_utack_restrict : forall a,
axioms a -> length a > 0 ->
a = utack (component (length a - 1) a) (restrict a (length a -1)).

Definition uple_map (f:E ->E) u :=
Uple.create (fun i => f (component i u)) (length u).

Lemma length_uple_map : forall f u,
length (uple_map f u) = length u.

Lemma component_uple_map : forall f u i,
i < length u -> component i (uple_map f u) = f (component i u).

Lemma axioms_uple_map : forall f u, axioms (uple_map f u).

Lemma uple_map_uple1 : forall f u,
uple_map f (uple1 u) = uple1 (f u).

Lemma uple_map_emptyset : forall f,
uple_map f emptyset = emptyset.

Lemma uple_map_concatenate : forall f u v,
uple_map f (concatenate u v) = concatenate (uple_map f u)
(uple_map f v).

End Uple.

\end{verbatim}

\subsection{The module {\tt Graph}}

\begin{verbatim}
Module Graph.
Export Uple.

Definition Vertices := R (v_(r_(t_ DOT ))).
Definition Edges := R (e_(d_(g_ DOT))).

Definition vertices a := V Vertices a.
Definition edges a := V Edges a.

Definition create v e :=
denote Vertices v
(denote Edges e stop).

Definition like a := a = create (vertices a) (edges a).

Lemma vertices_create : forall v e, vertices (create v e) =v.

Lemma edges_create : forall v e, edges (create v e) = e.

Lemma create_like : forall v e, like (create v e).

Lemma like_extensionality : forall a b,
like a -> like b -> vertices a = vertices b ->
edges a = edges b -> a = b.

Definition axioms a :=
like a &
(forall u, inc u (edges a) -> Arrow.like u) &
(forall u, inc u (edges a) -> inc (source u) (vertices a)) &
(forall u, inc u (edges a) -> inc (target u) (vertices a)).

Lemma axioms_extensionality : forall a b,
axioms a -> axioms b -> vertices a = vertices b ->
edges a = edges b -> a = b.

End Graph.

\end{verbatim}

\subsection{The module {\tt Free\_Category}}

\begin{verbatim}
Module Free_Category.
Export Graph.

Definition segment (i:nat) u :=
component i (arrow u).

Definition seg_length u :=
length (arrow u).

Definition arrow_chain u :=
Arrow.like u &
Uple.axioms (arrow u) &
(0 < seg_length u -> source (segment 0 u) = source u) &
(0 < seg_length u ->
target (segment (seg_length u -1) u) = target u) &
(forall i, i+1 < seg_length u ->
source (segment (i+1) u) =
target (segment i u)) &
(seg_length u = 0 -> source u = target u).

Definition mor_freecat a u :=
axioms a &
inc (source u) (vertices a) &
inc (target u) (vertices a) &
arrow_chain u &
(forall i, i < seg_length u ->
inc (segment i u) (edges a)).

Definition vertex_uple u :=
concatenate
(Uple.create (fun i => source (segment i u)) (seg_length u))
(uple1 (target u)).

Definition vertex_uple' u :=
concatenate (uple1 (source u))
(Uple.create (fun i => target (segment i u)) (seg_length u)).

Lemma length_vertex_uple : forall u,
length (vertex_uple u) = seg_length u + 1.

Lemma length_vertex_uple' : forall u,
length (vertex_uple' u) = seg_length u + 1.

Lemma vertex_uples_same : forall u,
arrow_chain u ->
vertex_uple u = vertex_uple' u.

Definition freecat_comp u v :=
Arrow.create (source v) (target u)
(concatenate (arrow v) (arrow u)).

Lemma seg_length_freecat_comp : forall u v,
seg_length (freecat_comp u v) = seg_length u + seg_length v.

Lemma segment_freecat_comp_first : forall u v i,
i < seg_length v -> segment i (freecat_comp u v)
= segment i v.

Lemma segment_freecat_comp_second : forall u v i,
seg_length v <= i -> i < seg_length u + seg_length v ->
segment i (freecat_comp u v) =
segment (i -seg_length v) u.
 

Definition freecat_id x :=
Arrow.create x x emptyset.

Lemma seg_length_freecat_id : forall x,
seg_length (freecat_id x) = 0.

Lemma arrow_chain_extensionality : forall u v,
arrow_chain u -> arrow_chain v ->
seg_length u = seg_length v ->
source u = source v -> target u = target v ->
(forall i, i< seg_length u -> segment i u = segment i v) ->
u = v.

Definition freecat_edge b :=
Arrow.create (source b) (target b) (uple1 b).

Lemma seg_length_freecat_edge : forall b,
seg_length (freecat_edge b) = 1.

Lemma segment_freecat_edge : forall b i,
i = 0 -> segment i (freecat_edge b) = b.

Lemma source_first_segment : forall u i,
arrow_chain u -> 0 < seg_length u -> i = 0 ->
source (segment i u) = source u.

Lemma target_last_segment : forall u i,
arrow_chain u -> 0 < seg_length u ->
i = seg_length u - 1 ->
target (segment i u) = target u.

Lemma eq_freecat_edge : forall u x,
arrow_chain u -> seg_length u = 1 ->
x = (segment 0 u) ->
u = freecat_edge x.

Lemma source_freecat_id : forall x,
source (freecat_id x) = x.

Lemma target_freecat_id : forall x,
target (freecat_id x) = x.

Lemma source_freecat_edge : forall b,
source (freecat_edge b) = source b.

Lemma target_freecat_edge : forall b,
target (freecat_edge b) = target b.

Lemma source_freecat_comp : forall u v,
source (freecat_comp u v) = source v.

Lemma target_freecat_comp : forall u v,
target (freecat_comp u v) = target u.

Lemma arrow_freecat_id : forall x,
arrow (freecat_id x) = emptyset.

Lemma arrow_freecat_edge : forall b,
arrow (freecat_edge b) = uple1 b.

Lemma arrow_freecat_comp : forall u v,
arrow (freecat_comp u v) =
concatenate (arrow v) (arrow u).

Lemma arrow_chain_freecat_id : forall x,
arrow_chain (freecat_id x).

Lemma eq_freecat_id : forall u,
arrow_chain u -> seg_length u = 0 ->
u = freecat_id (source u).

Lemma arrow_chain_freecat_edge : forall b,
arrow_chain (freecat_edge b).

Lemma arrow_chain_freecat_comp : forall u v,
arrow_chain u -> arrow_chain v -> source u = target v ->
arrow_chain (freecat_comp u v).

Lemma freecat_assoc : forall u v w,
freecat_comp u (freecat_comp v w) =
freecat_comp (freecat_comp u v) w.

Lemma freecat_left_id' : forall u x,
Arrow.like u -> x = target u -> Uple.axioms (arrow u) ->
freecat_comp (freecat_id x) u = u.

Lemma freecat_right_id' : forall u x,
Arrow.like u -> x = source u -> Uple.axioms (arrow u) ->
freecat_comp u (freecat_id x) = u.

Lemma freecat_left_id : forall u x,
arrow_chain u -> x = target u ->
freecat_comp (freecat_id x) u = u.

Lemma freecat_right_id : forall u x,
arrow_chain u -> x = source u ->
freecat_comp u (freecat_id x) = u.

Lemma uples_contained : forall b u,
Uple.axioms u ->
(forall i, i<length u -> inc (component i u) b) ->
inc u (powerset (Cartesian.product nat b)).

Lemma mor_freecat_from_set : forall a u,
mor_freecat a u ->
inc u (Image.create
(product (product (vertices a) (vertices a))
(powerset (product nat (edges a))))
(fun x => Arrow.create (pr1 (pr1 x)) (pr2 (pr1 x)) (pr2 x))).

Lemma mor_freecat_bounded : forall a,
Bounded.axioms (mor_freecat a).

Definition freecat_morphisms a :=
Bounded.create (mor_freecat a).

Lemma inc_freecat_morphisms : forall a u,
inc u (freecat_morphisms a) = mor_freecat a u.

Definition freecat a :=
Category.Notations.create
(vertices a) (freecat_morphisms a) freecat_comp freecat_id
emptyset.

Lemma is_mor_freecat : forall a u,
is_mor (freecat a) u = mor_freecat a u.

Lemma is_ob_freecat : forall a x,
is_ob (freecat a) x = inc x (vertices a).

Lemma comp_freecat1 : forall a u v,
mor_freecat a u -> mor_freecat a v ->
source u = target v ->
comp (freecat a) u v = freecat_comp u v.

Lemma id_freecat1 : forall a x,
inc x (vertices a) ->
id (freecat a) x = freecat_id x.

Lemma structure_freecat : forall a,
structure (freecat a) = emptyset.

Lemma mor_freecat_id : forall a x,
inc x (vertices a) -> axioms a ->
mor_freecat a (freecat_id x).

Lemma mor_freecat_comp : forall a u v,
mor_freecat a u -> mor_freecat a v ->
source u = target v ->
mor_freecat a (freecat_comp u v).

Lemma mor_freecat_edge : forall a u,
axioms a -> inc u (edges a) ->
mor_freecat a (freecat_edge u).

Lemma freecat_axioms : forall a,
axioms a -> Category.axioms (freecat a).

Lemma mor_freecat_rw : forall a u,
Graph.axioms a ->
mor (freecat a) u = mor_freecat a u.

Lemma ob_freecat_rw : forall a x,
Graph.axioms a ->
ob (freecat a) x = inc x (vertices a).

Lemma id_freecat : forall a x,
Graph.axioms a -> ob (freecat a) x ->
id (freecat a) x = freecat_id x.

Lemma comp_freecat : forall a u v,
Graph.axioms a -> mor (freecat a) u ->
mor (freecat a) v -> source u = target v ->
comp (freecat a) u v = freecat_comp u v.

Definition mor_chain a u :=
arrow_chain u &
Category.axioms a &
ob a (source u) & ob a (target u) &
(forall i, i< (seg_length u) -> mor a (segment i u)).

Lemma mor_chain_freecat_id : forall a x,
ob a x -> mor_chain a (freecat_id x).

Lemma mor_chain_freecat_comp : forall a u v,
mor_chain a u -> mor_chain a v ->
source u = target v -> mor_chain a (freecat_comp u v).

Lemma mor_chain_freecat_edge : forall a u,
mor a u -> mor_chain a (freecat_edge u).

Definition chain_tack u v :=
freecat_comp (freecat_edge u) v.

Lemma arrow_chain_chain_tack : forall u v,
arrow_chain v ->
source u = target v -> arrow_chain (chain_tack u v).

Lemma source_chain_tack : forall u v,
source (chain_tack u v) = source v.

Lemma target_chain_tack : forall u v,
target (chain_tack u v) = target u.

Lemma arrow_chain_tack : forall u v,
arrow (chain_tack u v) = utack u (arrow v).

Lemma seg_length_chain_tack : forall u v,
seg_length (chain_tack u v) = seg_length v + 1.

Lemma segment_chain_tack_old : forall u v i,
i < seg_length v ->
segment i (chain_tack u v) = segment i v.

Lemma segment_chain_tack_new : forall u v i,
i = seg_length v ->
segment i (chain_tack u v) = u.

Lemma mor_chain_chain_tack : forall a u v,
mor_chain a v -> mor a u ->
source u = target v -> mor_chain a (chain_tack u v).

Definition object_number i u :=
Y (i< seg_length u) (source (segment i u)) (target u).

Lemma object_number_zero : forall u,
arrow_chain u -> object_number 0 u = source u.

Lemma object_number_seg_length : forall u ,
arrow_chain u -> object_number (seg_length u) u = target u.

Lemma source_segment : forall i u,
arrow_chain u -> i < seg_length u ->
source (segment i u) = object_number i u.

Lemma target_segment : forall i u,
arrow_chain u -> i < seg_length u ->
target (segment i u) = object_number (i+1) u.

Definition chain_restrict i u :=
Arrow.create (source u) (object_number i u)
(restrict (arrow u) i).

Lemma source_chain_restrict : forall i u,
source (chain_restrict i u) = source u.

Lemma target_chain_restrict : forall i u,
target (chain_restrict i u) = object_number i u.

Lemma seg_length_chain_restrict : forall i u,
seg_length (chain_restrict i u) = i.

Lemma segment_chain_restrict : forall i j u,
j < i ->
segment j (chain_restrict i u) = segment j u.

Lemma arrow_chain_chain_restrict : forall i u,
arrow_chain u -> i <= seg_length u ->
arrow_chain (chain_restrict i u).

Lemma chain_restrict_refl : forall u,
arrow_chain u ->
chain_restrict (seg_length u) u = u.

Lemma eq_chain_tack_restrict : forall u,
arrow_chain u -> seg_length u > 0 ->
u = chain_tack (segment (seg_length u - 1) u)
(chain_restrict (seg_length u - 1) u).

Lemma arrow_chain_induction : forall (P:E->Prop),
(forall x, P (freecat_id x)) ->
(forall u v, arrow_chain v -> source u = target v ->
P v -> P (chain_tack u v)) ->
(forall y, arrow_chain y -> P y).

Lemma mor_chain_chain_restrict : forall a u i,
mor_chain a u -> i <= seg_length u ->
mor_chain a (chain_restrict i u).

Lemma chain_restrict_chain_tack : forall u v i,
i = seg_length v -> arrow_chain v ->
source u = target v ->
chain_restrict i (chain_tack u v) = v.

Lemma mor_chain_induction : forall a (P:E->Prop),
(forall x, ob a x -> P (freecat_id x)) ->
(forall u v, mor_chain a v -> mor a u ->
source u = target v -> P v -> P (chain_tack u v)) ->
(forall y, mor_chain a y -> P y).

Definition compose_chain a u :=
iterate (seg_length u) (fun i v => comp a (segment i u) v)
(id a (source u)).

Lemma compose_chain_freecat_id : forall a x,
compose_chain a (freecat_id x) = id a x.

Lemma compose_chain_chain_tack : forall a u v,
compose_chain a (chain_tack u v) =
comp a u (compose_chain a v).

Lemma compose_chain_freecat_edge : forall a u,
mor a u -> compose_chain a (freecat_edge u) = u.

Definition compose_chain_facts a u :=
mor_chain a u &
source (compose_chain a u) = source u &
target (compose_chain a u) = target u &
mor a (compose_chain a u).

Lemma compose_chain_facts_compose_chain : forall a u,
mor_chain a u -> compose_chain_facts a u.

Lemma source_compose_chain : forall a u,
mor_chain a u -> source (compose_chain a u) = source u.

Lemma target_compose_chain : forall a u,
mor_chain a u -> target (compose_chain a u) = target u.

Lemma mor_compose_chain : forall a u, mor_chain a u ->
mor a (compose_chain a u).

Lemma freecat_comp_chain_tack : forall y u v,
freecat_comp (chain_tack y u) v =
chain_tack y (freecat_comp u v).

Lemma compose_chain_freecat_comp : forall a u v,
mor_chain a u -> mor_chain a v ->
source u = target v ->
compose_chain a (freecat_comp u v) =
comp a (compose_chain a u) (compose_chain a v).

Definition chain_map (fo : E -> E) (fa:E -> E) u :=
Arrow.create
(fo (source u)) (fo (target u)) (uple_map fa (arrow u)).

Lemma source_chain_map : forall fo fa u,
source (chain_map fo fa u) = fo (source u).

Lemma target_chain_map : forall fo fa u,
target (chain_map fo fa u) = fo (target u).

Lemma arrow_chain_map : forall fo fa u,
arrow (chain_map fo fa u) = uple_map fa (arrow u).

Lemma seg_length_chain_map : forall fo fa u,
seg_length (chain_map fo fa u) = seg_length u.

Lemma segment_chain_map : forall fo fa u i,
i < seg_length u -> segment i (chain_map fo fa u)
= fa (segment i u).

Lemma arrow_chain_chain_map : forall fo fa u,
arrow_chain u ->
(forall i, i<seg_length u -> fo (source (segment i u)) =
source (fa (segment i u))) ->
(forall i, i<seg_length u -> fo (target (segment i u)) =
target (fa (segment i u))) ->
arrow_chain (chain_map fo fa u).

Lemma chain_map_freecat_id : forall fo fm x,
chain_map fo fm (freecat_id x) = freecat_id (fo x).

Lemma chain_map_freecat_edge : forall fo fm u,
source (fm u) = fo (source u) ->
target (fm u) = fo (target u) ->
chain_map fo fm (freecat_edge u) = freecat_edge (fm u).

Lemma chain_map_freecat_comp : forall fo fm u v,
chain_map fo fm (freecat_comp u v) =
freecat_comp (chain_map fo fm u) (chain_map fo fm v).

Lemma chain_map_chain_tack : forall fo fm u v,
target (fm u) = fo (target u) ->
chain_map fo fm (chain_tack u v) =
chain_tack (fm u) (chain_map fo fm v).

Definition free_functor g a fo fm :=
Functor.create (freecat g) a (fun u =>
(compose_chain a (chain_map fo fm u))).

Definition free_functor_property g a fo fm :=
Graph.axioms g &
Category.axioms a &
(forall x, inc x (vertices g) -> ob a (fo x)) &
(forall u, inc u (edges g) -> mor a (fm u)) &
(forall u, inc u (edges g) -> source (fm u) = fo (source u)) &
(forall u, inc u (edges g) -> target (fm u) = fo (target u)).

Lemma source_free_functor : forall g a fo fm,
source (free_functor g a fo fm) = freecat g.

Lemma target_free_functor : forall g a fo fm,
target (free_functor g a fo fm) = a.

Lemma fmor_free_functor : forall g a fo fm u,
Graph.axioms g ->
mor_freecat g u -> fmor (free_functor g a fo fm) u =
(compose_chain a (chain_map fo fm u)).

Lemma mor_chain_chain_map : forall g a fo fm u,
free_functor_property g a fo fm ->
mor_freecat g u ->
mor_chain a (chain_map fo fm u).

Lemma fob_free_functor : forall g a fo fm x,
free_functor_property g a fo fm ->
inc x (vertices g) ->
fob (free_functor g a fo fm) x = fo x.

Lemma free_functor_axioms : forall g a fo fm,
free_functor_property g a fo fm ->
Functor.axioms (free_functor g a fo fm).

Lemma fmor_ff_freecat_id : forall g a fo fm x,
free_functor_property g a fo fm ->
inc x (vertices g) ->
fmor (free_functor g a fo fm) (freecat_id x) = id a (fo x).

Lemma fmor_ff_freecat_edge : forall g a fo fm u,
free_functor_property g a fo fm ->
inc u (edges g) ->
fmor (free_functor g a fo fm) (freecat_edge u) = fm u.

Lemma fmor_comp : forall a f u v,
Functor.axioms f -> source f = a ->
mor a u -> mor a v ->
source u = target v ->
fmor f (comp a u v) = comp (target f) (fmor f u) (fmor f v).

Lemma fmor_ff_freecat_comp : forall g a fo fm u v,
free_functor_property g a fo fm ->
mor (freecat g) u ->
mor (freecat g) v ->
source u = target v ->
fmor (free_functor g a fo fm) (freecat_comp u v) =
comp a (fmor (free_functor g a fo fm) u)
(fmor (free_functor g a fo fm) v).

Lemma fmor_ff_chain_tack : forall g a fo fm u v,
free_functor_property g a fo fm ->
inc u (edges g) ->
mor (freecat g) v -> source u = target v ->
fmor (free_functor g a fo fm) (chain_tack u v) =
comp a (fm u) (fmor (free_functor g a fo fm) v).

Lemma free_functor_property_fob_fmor : forall g a f,
Functor.axioms f -> Graph.axioms g ->
source f = freecat g -> target f = a ->
free_functor_property g a (fob f) 
(fun u => fmor f (freecat_edge u)).

Definition graph_fmor g f y :=
fmor (free_functor g (target f) (fob f)
(fun u => fmor f (freecat_edge u))) y.

Lemma graph_fob : forall g f x,
Graph.axioms g -> source f = freecat g ->
Functor.axioms f -> inc x (vertices g) ->
fob (free_functor g (target f) (fob f)
(fun u => fmor f (freecat_edge u))) x = fob f x.

Definition graph_fmor_recovers g f y :=
graph_fmor g f y = fmor f y.

Lemma graph_fmor_recovers_freecat_id : forall g f x,
Graph.axioms g -> source f = freecat g ->
Functor.axioms f -> inc x (vertices g) ->
graph_fmor_recovers g f (freecat_id x).

Lemma graph_fmor_recovers_freecat_edge : forall g f u,
Graph.axioms g -> source f = freecat g ->
Functor.axioms f -> inc u (edges g) ->
graph_fmor_recovers g f (freecat_edge u).

Lemma graph_fmor_recovers_freecat_comp : forall g f u v,
Graph.axioms g -> source f = freecat g ->
Functor.axioms f -> mor_freecat g u -> mor_freecat g v ->
source u = target v ->
graph_fmor_recovers g f u -> graph_fmor_recovers g f v ->
graph_fmor_recovers g f (freecat_comp u v).

Lemma graph_fmor_recovers_chain_tack : forall g f u v,
Graph.axioms g -> source f = freecat g ->
Functor.axioms f -> inc u (edges g) -> mor_freecat g v ->
source u = target v ->
graph_fmor_recovers g f v ->
graph_fmor_recovers g f (chain_tack u v).

Lemma inc_object_number_vertices : forall g u i,
Graph.axioms g -> mor_freecat g u ->
i <= seg_length u -> inc (object_number i u) (vertices g).

Lemma mor_freecat_chain_restrict : forall g u i,
Graph.axioms g -> mor_freecat g u ->
i <= seg_length u -> mor_freecat g (chain_restrict i u).

Lemma mor_freecat_induction : forall g (P : E -> Prop),
Graph.axioms g ->
(forall x, inc x (vertices g) -> P (freecat_id x)) ->
(forall u v, inc u (edges g) -> mor_freecat g v ->
source u = target v -> P v -> P (chain_tack u v)) ->
(forall u, mor_freecat g u -> P u).

Lemma graph_fmor_recovers_all : forall g f u,
Graph.axioms g -> source f = freecat g ->
Functor.axioms f -> mor_freecat g u ->
graph_fmor_recovers g f u.

Lemma eq_free_functor_fob_fmor : forall g f,
Graph.axioms g -> source f = freecat g ->
Functor.axioms f ->
f = free_functor g (target f) (fob f)
(fun u => fmor f (freecat_edge u)).

Definition free_nt_property g a b (t:E -> E) :=
Graph.axioms g & Functor.axioms a & Functor.axioms b
& source a = freecat g
& source b = freecat g
& target a = target b &
(forall x, inc x (vertices g) -> mor (target a) (t x)) &
(forall x, inc x (vertices g) -> source (t x) = fob a x) &
(forall x, inc x (vertices g) -> target (t x) = fob b x) &
(forall u, inc u (edges g) ->
comp (target a) (t (target u)) (fmor a (freecat_edge u)) =
comp (target a) (fmor b (freecat_edge u)) (t (source u))).

Definition free_nt_respects a b (t:E->E) u :=
comp (target a) (t (target u)) (fmor a u) =
comp (target a) (fmor b u) (t (source u)).

Lemma free_nt_respects_freecat_id : forall g a b t x,
free_nt_property g a b t ->
inc x (vertices g) ->
free_nt_respects a b t (freecat_id x).

Lemma free_nt_respects_freecat_edge : forall g a b t u,
free_nt_property g a b t ->
inc u (edges g) ->
free_nt_respects a b t (freecat_edge u).

Lemma free_nt_respects_freecat_comp : forall g a b t u v,
free_nt_property g a b t ->
mor_freecat g u -> mor_freecat g v -> source u = target v ->
free_nt_respects a b t u ->
free_nt_respects a b t v ->
free_nt_respects a b t (freecat_comp u v).

Lemma free_nt_respects_chain_tack : forall g a b t u v,
free_nt_property g a b t -> inc u (edges g) ->
mor_freecat g v -> source u = target v ->
free_nt_respects a b t v ->
free_nt_respects a b t (chain_tack u v).

Lemma free_nt_respects_all : forall g a b t u,
free_nt_property g a b t ->
mor_freecat g u ->
free_nt_respects a b t u.

Lemma free_nt_nat_trans_property : forall g a b t,
free_nt_property g a b t -> Nat_Trans.property a b t.

Lemma free_nt_prop_eq_nat_trans_prop : forall g a b t,
Graph.axioms g -> source a = freecat g ->
free_nt_property g a b t = Nat_Trans.property a b t.

Lemma free_nt_property_ntrans : forall g r,
Graph.axioms g -> Nat_Trans.axioms r ->
osource r = freecat g ->
free_nt_property g (source r) (target r) (ntrans r).


End Free_Category.

\end{verbatim}

\newpage

\section{The file {\tt qcat.v}}

\begin{verbatim}

Require Export updateA.

\end{verbatim}

\subsection{The module {\tt Categorical\_Relation}}

\begin{verbatim}
Module Categorical_Relation.
Export UpdateA.
Export Category.

Definition cat_rel a r :=
Category.axioms a &
is_relation r &
(forall x y, related r x y -> mor a x) &
(forall x y, related r x y -> mor a y) &
(forall x y, related r x y -> source x = source y) &
(forall x y, related r x y -> target x = target y).

Definition cat_equiv_rel a r :=
cat_rel a r &
is_equivalence_relation r &
(forall u, mor a u -> related r u u) &
(forall x y u v, related r x y -> related r u v ->
source x = target u ->
related r (comp a x u) (comp a y v)).

Lemma cat_rel_related_rw : forall a r x y,
cat_rel a r ->
related r x y =
(related r x y & mor a x & mor a y &
source x = source y & target x = target y).

Lemma cer_reflexive : forall a r u,
cat_equiv_rel a r -> mor a u -> related r u u.

Lemma related_comp : forall a r x y u v,
cat_equiv_rel a r ->
related r x y -> related r u v -> source x = target u ->
related r (comp a x u) (comp a y v).

Lemma cat_rel_intersection : forall a z,
nonempty z ->
(forall r, inc r z -> cat_rel a r) ->
cat_rel a (intersection z).

Lemma cer_intersection : forall a z,
nonempty z ->
(forall r, inc r z -> cat_equiv_rel a r) ->
cat_equiv_rel a (intersection z).

Definition coarse a :=
Z (Cartesian.product (morphisms a) (morphisms a))
(fun u => (source (pr1 u) = source (pr2 u) &
target (pr1 u) = target (pr2 u))).

Lemma related_coarse : forall a x y,
Category.axioms a ->
related (coarse a) x y = (mor a x & mor a y &
source x = source y & target x = target y).

Lemma cat_rel_coarse : forall a,
Category.axioms a ->
cat_rel a (coarse a).

Lemma cat_equiv_rel_coarse : forall a,
Category.axioms a ->
cat_equiv_rel a (coarse a).

Lemma sub_coarse : forall a r,
cat_rel a r ->
sub r (coarse a).

Definition relations_between r s :=
Z (powerset s) (fun y => sub r y).

Lemma inc_relations_between : forall r s y,
is_relation s ->
inc y (relations_between r s) =
(sub r y & sub y s & is_relation y).

Definition cer a r :=
intersection (Z (relations_between r (coarse a))
(fun y=> cat_equiv_rel a y)).

Lemma cat_equiv_rel_cer : forall a r,
cat_rel a r -> cat_equiv_rel a (cer a r).

Lemma cat_rel_subset : forall a r,
(exists s, (cat_rel a s & sub r s)) ->
cat_rel a r.

Lemma sub_cer : forall a r s,
sub r s -> cat_equiv_rel a s ->
sub (cer a r) s.

Lemma cer_contains : forall a r,
cat_rel a r -> sub r (cer a r).

Definition ker f :=
Z (coarse (source f)) (fun z => (fmor f (pr1 z) = fmor f (pr2 z))).

Lemma related_ker : forall f x y,
Functor.axioms f ->
related (ker f) x y =
(mor (source f) x & mor (source f) y &
source x = source y & target x = target y & fmor f x = fmor f y).

Lemma is_relation_ker : forall f, Functor.axioms f ->
is_relation (ker f).

Lemma cat_rel_ker : forall a f,
Functor.axioms f -> a = source f ->
cat_rel a (ker f).

Lemma cat_equiv_rel_ker : forall a f,
Functor.axioms f -> a = source f ->
cat_equiv_rel a (ker f).

Definition compatible r f :=
Functor.axioms f & cat_rel (source f) r &
sub r (ker f).

Lemma compatible_related_ker : forall r f u v,
compatible r f -> related r u v -> related (ker f) u v.

Lemma related_eq_fmor : forall r f u v,
compatible r f -> related r u v -> fmor f u = fmor f v.

Lemma is_relation_coarse : forall a,
Category.axioms a -> is_relation (coarse a).

Lemma compatible_cer : forall a r f,
compatible r f -> a = source f ->
compatible (cer a r) f.

Lemma compatible_rw : forall r f,
compatible r f =
(Functor.axioms f & cat_rel (source f) r &
(forall x y, related r x y -> fmor f x = fmor f y)).

End Categorical_Relation.

\end{verbatim}

\subsection{The module {\tt Quotient\_Category}}

\begin{verbatim}
Module Quotient_Category.
Export Categorical_Relation.

Definition arrow_class r u := 
Arrow.create (source u) (target u) (class r u).

Lemma source_arrow_class : forall r u,
source (arrow_class r u) = source u.

Lemma target_arrow_class : forall r u,
target (arrow_class r u) = target u.

Lemma arrow_arrow_class : forall r u,
arrow (arrow_class r u) = class r u.

Lemma inc_arrow_arrow_class : forall a r u v,
cat_equiv_rel a r -> mor a v ->
inc v (arrow (arrow_class r u)) = related r u v.

Lemma like_arrow_class : forall r u,
Arrow.like (arrow_class r u).

Definition is_quotient_arrow a r u :=
cat_equiv_rel a r &
(exists y, mor a y & u = arrow_class r y).

Lemma is_quotient_arrow_arrow_class : forall a r u,
cat_equiv_rel a r -> mor a u ->
is_quotient_arrow a r (arrow_class r u).

Definition quotient_morphisms a r :=
Image.create (morphisms a) (arrow_class r).

Lemma inc_quotient_morphisms : forall a r u,
cat_equiv_rel a r ->
inc u (quotient_morphisms a r) = is_quotient_arrow a r u.

Lemma related_arrow_class_eq : forall a r u v,
cat_equiv_rel a r ->
mor a u ->
related r u v =
(arrow_class r u = arrow_class r v).

Definition arrow_rep v := rep (arrow v).

Lemma inc_arrow_class_refl : forall a r u,
cat_equiv_rel a r -> mor a u -> inc u (arrow (arrow_class r u)).

Lemma nonempty_arrow : forall a r u,
is_quotient_arrow a r u ->
nonempty (arrow u).

Lemma inc_arrow_rep_arrow : forall a r u,
is_quotient_arrow a r u ->
inc (arrow_rep u) (arrow u).

Lemma related_arrow_rep_arrow_class : forall a r u,
cat_equiv_rel a r -> mor a u ->
related r u (arrow_rep (arrow_class r u)).

Lemma source_arrow_rep : forall a r u,
is_quotient_arrow a r u -> source (arrow_rep u) = source u.

Lemma target_arrow_rep : forall a r u,
is_quotient_arrow a r u -> target (arrow_rep u) = target u.

Lemma arrow_class_arrow_rep : forall a r u,
is_quotient_arrow a r u ->
arrow_class r (arrow_rep u) = u.

Lemma inc_arrow_facts : forall a r u y,
is_quotient_arrow a r u -> inc y (arrow u) ->
(mor a y & source y = source u & target y = target u).

Lemma mor_arrow_rep : forall a r u,
is_quotient_arrow a r u ->
mor a (arrow_rep u).

Lemma related_arrow_rep : forall a r u v,
cat_equiv_rel a r -> mor a v -> u = arrow_class r v ->
related r v (arrow_rep u).

Lemma related_arrow_rep_rw : forall a r u v,
cat_equiv_rel a r -> is_quotient_arrow a r u -> mor a v ->
related r v (arrow_rep u) = (u = arrow_class r v).

Definition quot_id a r x :=
arrow_class r (id a x).

Definition quot_comp a r u v :=
arrow_class r (comp a (arrow_rep u) (arrow_rep v)).

Lemma source_quot_id : forall a r x,
ob a x ->
source (quot_id a r x) = x.

Lemma target_quot_id : forall a r x,
ob a x ->
target (quot_id a r x) = x.

Lemma source_quot_comp : forall a r u v,
cat_equiv_rel a r ->
is_quotient_arrow a r u ->
is_quotient_arrow a r v -> source u = target v ->
source (quot_comp a r u v) = source v.

Lemma target_quot_comp : forall a r u v,
cat_equiv_rel a r ->
is_quotient_arrow a r u ->
is_quotient_arrow a r v -> source u = target v ->
target (quot_comp a r u v) = target u.

Lemma is_quotient_arrow_quot_id : forall a r x,
cat_equiv_rel a r -> ob a x ->
is_quotient_arrow a r (quot_id a r x).

Lemma is_quotient_arrow_quot_comp : forall a r u v,
cat_equiv_rel a r -> is_quotient_arrow a r u -> 
is_quotient_arrow a r v -> source u = target v ->
is_quotient_arrow a r (quot_comp a r u v).

Lemma arrow_class_eq : forall a r u v,
cat_equiv_rel a r -> mor a u ->
(arrow_class r u = arrow_class r v) =
(related r u v).

Lemma arrow_class_id : forall a r x,
arrow_class r (id a x) = quot_id a r x.

Lemma quot_comp_arrow_class : forall a r u v,
cat_equiv_rel a r -> mor a u -> mor a v -> source u = target v ->
quot_comp a r (arrow_class r u) (arrow_class r v) =
arrow_class r (comp a u v).

Lemma quot_id_left : forall a r a' r' x u,
cat_equiv_rel a r -> is_quotient_arrow a r u ->
x = target u -> a' = a -> r' = r ->
quot_comp a r (quot_id a' r' x) u = u.

Lemma quot_id_right : forall a r a' r' x u,
cat_equiv_rel a r -> is_quotient_arrow a r u ->
x = source u -> a' = a -> r' = r ->
quot_comp a r u (quot_id a' r' x) = u.

Lemma quot_comp_assoc : forall a r a' r' u v w,
cat_equiv_rel a r -> is_quotient_arrow a r u -> 
is_quotient_arrow a r v ->
is_quotient_arrow a r w -> source u = target v ->
source v = target w -> a' = a -> r' = r ->
quot_comp a r (quot_comp a' r' u v) w = 
quot_comp a r u (quot_comp a' r' v w).

Definition quotient_cat a r := Category.Notations.create (objects a)
(quotient_morphisms a r) (quot_comp a r) (quot_id a r) (structure a).

Lemma is_ob_quotient_cat : forall a r x,
is_ob (quotient_cat a r) x = is_ob a x.

Lemma is_mor_quotient_cat : forall a r u,
cat_equiv_rel a r ->
is_mor (quotient_cat a r) u = is_quotient_arrow a r u.

Lemma comp_quotient_cat : forall a r u v,
cat_equiv_rel a r -> is_quotient_arrow a r u -> 
is_quotient_arrow a r v ->
source u = target v -> 
comp (quotient_cat a r) u v = quot_comp a r u v.

Lemma id_quotient_cat : forall a r x,
cat_equiv_rel a r -> ob a x ->
id (quotient_cat a r) x = quot_id a r x.

Lemma quotient_cat_axioms : forall a r,
cat_equiv_rel a r ->
Category.axioms (quotient_cat a r).

Lemma ob_quotient_cat : forall a r x,
cat_equiv_rel a r ->
ob (quotient_cat a r) x = ob a x.

Lemma mor_quotient_cat : forall a r u,
cat_equiv_rel a r ->
mor (quotient_cat a r) u = is_quotient_arrow a r u.

Lemma mor_quotient_cat_ex : forall a r u,
cat_equiv_rel a r ->
mor (quotient_cat a r) u =
(exists y, mor a y & u = arrow_class r y).

Lemma mor_quotient_cat_quot_id : forall a r x,
cat_equiv_rel a r -> ob a x ->
mor (quotient_cat a r) (quot_id a r x).

Lemma mor_quotient_cat_quot_comp : forall a r u v,
cat_equiv_rel a r ->
mor (quotient_cat a r) u ->
mor (quotient_cat a r) v ->
source u = target v ->
mor (quotient_cat a r) (quot_comp a r u v).

Lemma mor_quotient_cat_arrow_class : forall a r u,
cat_equiv_rel a r -> mor a u ->
mor (quotient_cat a r) (arrow_class r u).

End Quotient_Category.

\end{verbatim}

\subsection{The module {\tt Quotient\_Functor}}

\begin{verbatim}
Module Quotient_Functor.
Export Quotient_Category.

Definition qfunctor a b r (fo fm:E-> E) :=
Functor.create a (quotient_cat b r) (fun u => arrow_class r (fm u)).

Lemma source_qfunctor : forall a b r fo fm,
source (qfunctor a b r fo fm) = a.

Lemma target_qfunctor : forall a b r fo fm,
target (qfunctor a b r fo fm) = (quotient_cat b r).

Lemma fmor_qfunctor : forall a b r fo fm u,
mor a u ->
fmor (qfunctor a b r fo fm) u = arrow_class r (fm u).

Definition qfunctor_property a b r fo fm :=
Category.axioms a &
cat_equiv_rel b r &
(forall x, ob a x -> ob b (fo x)) &
(forall u, mor a u -> mor b (fm u)) &
(forall u, mor a u -> source (fm u) = fo (source u)) &
(forall u, mor a u -> target (fm u) = fo (target u)) &
(forall x, ob a x -> related r (fm (id a x)) (id b (fo x))) &
(forall u v, mor a u -> mor a v -> source u = target v ->
related r (fm (comp a u v)) (comp b (fm u) (fm v))).

Lemma fob_qfunctor : forall a b r fo fm x,
qfunctor_property a b r fo fm -> ob a x ->
fob (qfunctor a b r fo fm) x = fo x.

Lemma ob_fob_qfunctor : forall a b r fo fm x,
qfunctor_property a b r fo fm -> ob a x ->
ob (quotient_cat b r) (fob (qfunctor a b r fo fm) x).

Lemma mor_fmor_qfunctor : forall a b r fo fm u,
qfunctor_property a b r fo fm -> mor a u ->
mor (quotient_cat b r) (fmor (qfunctor a b r fo fm) u).

Lemma source_fmor_qfunctor : forall a b r fo fm u,
qfunctor_property a b r fo fm -> mor a u ->
source (fmor (qfunctor a b r fo fm) u) = 
fob (qfunctor a b r fo fm) (source u).

Lemma target_fmor_qfunctor : forall a b r fo fm u,
qfunctor_property a b r fo fm -> mor a u ->
target (fmor (qfunctor a b r fo fm) u) = 
fob (qfunctor a b r fo fm) (target u).

Lemma fmor_qfunctor_id : forall a b r fo fm x,
qfunctor_property a b r fo fm -> ob a x ->
fmor (qfunctor a b r fo fm) (id a x) =
id (quotient_cat b r) (fob (qfunctor a b r fo fm) x).

Lemma fmor_qfunctor_comp : forall a b r fo fm u v,
qfunctor_property a b r fo fm -> mor a u ->
mor a v -> source u = target v ->
fmor (qfunctor a b r fo fm) (comp a u v) =
comp (quotient_cat b r) (fmor (qfunctor a b r fo fm) u)
(fmor (qfunctor a b r fo fm) v).

Lemma qfunctor_axioms : forall a b r fo fm,
qfunctor_property a b r fo fm ->
Functor.axioms (qfunctor a b r fo fm).

Lemma qfunctor_extensionality : forall a b r fo fm a' b' r' fo' fm',
a = a' -> b = b' -> r = r' ->
qfunctor_property a b r fo fm ->
qfunctor_property a' b' r' fo' fm' ->
(forall u, mor a u -> related r (fm u) (fm' u)) ->
qfunctor a b r fo fm = qfunctor a' b' r' fo' fm'.

Lemma eq_qfunctor : forall a b r f fo fm,
Functor.axioms f -> source f = a ->
target f = (quotient_cat b r) ->
qfunctor_property a b r fo fm ->
(forall u, mor a u -> related r (fm u) (arrow_rep (fmor f u))) ->
f = qfunctor a b r fo fm.

Lemma qfunctor_property_fob_fmor : forall b r f,
Functor.axioms f ->
cat_equiv_rel b r -> target f = (quotient_cat b r) ->
qfunctor_property (source f) b r 
(fob f) (fun u => arrow_rep (fmor f u)).

Lemma eq_qfunctor2 : forall b r f,
Functor.axioms f ->
cat_equiv_rel b r -> target f = (quotient_cat b r) ->
f = qfunctor (source f) b r (fob f) (fun u => arrow_rep (fmor f u)).

Definition qprojection a r :=
qfunctor a a r (fun (x:E) => x) (fun (u:E)=> u).

Lemma source_qprojection : forall a r,
source (qprojection a r) = a.

Lemma target_qprojection : forall a r,
target (qprojection a r) = (quotient_cat a r).

Lemma qprojection_property : forall a r,
cat_equiv_rel a r ->
qfunctor_property a a r (fun (x:E) => x) (fun (u:E)=> u).

Lemma fob_qprojection : forall a r x,
cat_equiv_rel a r -> ob a x ->
fob (qprojection a r) x = x.

Lemma fmor_qprojection : forall a r u,
mor a u ->
fmor (qprojection a r) u = (arrow_class r u).

Lemma fmor_qprojection_arrow_rep : forall a r u,
is_quotient_arrow a r u ->
fmor (qprojection a r) (arrow_rep u) = u.

Lemma qprojection_axioms : forall a r,
cat_equiv_rel a r ->
Functor.axioms (qprojection a r).

Definition qdotted r f :=
Functor.create (quotient_cat (source f) r) (target f)
(fun u => fmor f (arrow_rep u)).

Lemma source_qdotted : forall r f,
source (qdotted r f) = quotient_cat (source f) r.

Lemma target_qdotted : forall r f,
target (qdotted r f) = target f.

Lemma fmor_qdotted : forall r f u,
cat_equiv_rel (source f) r
-> is_quotient_arrow (source f) r u ->
fmor (qdotted r f) u = fmor f (arrow_rep u).

Lemma fmor_qdotted_arrow_class : forall r f u,
cat_equiv_rel (source f) r ->
compatible r f -> mor (source f) u ->
fmor (qdotted r f) (arrow_class r u) = fmor f u.

Lemma fob_qdotted : forall r f x,
cat_equiv_rel (source f) r ->
compatible r f -> ob (source f) x ->
fob (qdotted r f) x = fob f x.

Lemma fmor_qdotted_quot_id : forall r f x,
cat_equiv_rel (source f) r ->
compatible r f -> ob (source f) x ->
fmor (qdotted r f) (quot_id (source f) r x) = id (target f) (fob f x).

Lemma fmor_qdotted_quot_comp : forall r f u v,
cat_equiv_rel (source f) r ->
compatible r f -> is_quotient_arrow (source f) r u ->
is_quotient_arrow (source f) r v -> source u = target v ->
fmor (qdotted r f) (quot_comp (source f) r u v) =
comp (target f) (fmor (qdotted r f) u) (fmor (qdotted r f) v).

Lemma qdotted_axioms : forall r f,
cat_equiv_rel (source f) r ->
compatible r f -> Functor.axioms (qdotted r f).

Lemma eq_fcompose_qdotted_qprojection : forall r f,
cat_equiv_rel (source f) r ->
compatible r f -> f = 
fcompose (qdotted r f) (qprojection (source f) r).

Lemma compatible_fcompose_qprojection : forall a r f,
Functor.axioms f -> cat_equiv_rel a r ->
source f = quotient_cat a r ->
compatible r (fcompose f (qprojection a r)).

Lemma eq_qdotted : forall a r f,
Functor.axioms f -> cat_equiv_rel a r ->
source f = quotient_cat a r ->
f = qdotted r (fcompose f (qprojection a r)).

End Quotient_Functor.

\end{verbatim}

\subsection{The module {\tt Ob\_Iso\_Functor}}

\begin{verbatim}
Module Ob_Iso_Functor.
Export Functor_Cat.
Export Quotient_Functor.

Definition pull_morphism a f :=
Functor.create (functor_cat (target f) a) 
(functor_cat (source f) a) (fun u => htrans_right u f).

Lemma source_pull_morphism : forall f a,
source (pull_morphism a f) = (functor_cat (target f) a).

Lemma target_pull_morphism : forall f a,
target (pull_morphism a f) = (functor_cat (source f) a).

Lemma fmor_pull_morphism : forall f a u,

mor (functor_cat (target f) a) u ->
fmor (pull_morphism a f) u = htrans_right u f.

Lemma fob_pull_morphism : forall f a y,
Functor.axioms f -> Category.axioms a ->
ob (functor_cat (target f) a) y ->
fob (pull_morphism a f) y = fcompose y f.

Lemma pull_morphism_axioms : forall f a,
Functor.axioms f -> Category.axioms a ->
Functor.axioms (pull_morphism a f).

Definition faithful f :=
Functor.axioms f &
(forall u v, related (ker f) u v -> u = v).

Definition says f x :=
Functor.axioms f &
(exists y, (ob (source f) y & fob f y = x)).

Definition msays f u :=
Functor.axioms f &
(exists y, (mor (source f) y & fmor f y = u)).

Definition full f :=
Functor.axioms f &
(forall u, mor (target f) u ->
says f (source u) -> says f (target u) -> msays f u).

Definition ob_inj f :=
Functor.axioms f &
(forall x y, ob (source f) x -> ob (source f) y -> 
fob f x = fob f y -> x = y).

Definition iso_to_full_subcategory f :=
faithful f & full f & ob_inj f.

Definition ob_image f :=
Z (objects (target f))
(fun x => (exists y, (ob (source f) y & fob f y = x))).

Definition ob_surj f :=
Functor.axioms f &
(forall x, ob (target f) x -> says f x).

Definition ob_iso f :=
ob_inj f & ob_surj f.

Lemma has_finverse_rw : forall f,
has_finverse f = (faithful f & full f & ob_iso f).

Definition is_full_subcategory a b :=
is_subcategory a b &
(forall u, mor b u -> ob a (source u) -> ob a (target u) ->
mor a u).

Definition full_subcategory (a:E) (obp : E -> Prop) :=
subcategory a obp (fun u => (obp (source u) & obp (target u))).

Lemma full_subcategory_property : forall a obp,
Category.axioms a -> subcategory_property a obp
(fun u => (obp (source u) & obp (target u))).

Lemma full_subcategory_axioms : forall a obp,
Category.axioms a ->
Category.axioms (full_subcategory a obp).

Lemma is_subcategory_full_subcategory : forall a obp,
Category.axioms a ->
is_subcategory (full_subcategory a obp) a.

Lemma ob_full_subcategory : forall a obp x,
axioms a ->
ob (full_subcategory a obp) x = (ob a x & obp x).

Lemma mor_full_subcategory : forall a obp u,
axioms a ->
mor (full_subcategory a obp) u = (mor a u & obp (source u)
& obp (target u)).

Lemma id_full_subcategory : forall a obp x,
axioms a ->
ob (full_subcategory a obp) x ->
id (full_subcategory a obp) x = id a x.

Lemma comp_full_subcategory : forall a obp u v,
axioms a ->
mor (full_subcategory a obp) u ->
mor (full_subcategory a obp) v ->
source u = target v ->
comp (full_subcategory a obp) u v = comp a u v.

Lemma is_full_subcategory_full_subcategory : forall a obp,
Category.axioms a ->
is_full_subcategory (full_subcategory a obp) a.

Definition subcategory_inclusion a b :=
Functor.create a b (fun u => u).

Lemma source_subcategory_inclusion : forall a b,
source (subcategory_inclusion a b) = a.

Lemma target_subcategory_inclusion : forall a b,
target (subcategory_inclusion a b) = b.

Lemma fmor_subcategory_inclusion : forall a b u,
is_subcategory a b -> mor a u ->
fmor (subcategory_inclusion a b) u = u.

Lemma fob_subcategory_inclusion : forall a b x,
is_subcategory a b -> ob a x ->
fob (subcategory_inclusion a b) x = x.

Lemma subcategory_inclusion_axioms : forall a b,
is_subcategory a b ->
Functor.axioms (subcategory_inclusion a b).

Lemma faithful_subcategory_inclusion : forall a b,
is_subcategory a b ->
faithful (subcategory_inclusion a b).

Lemma ob_inj_subcategory_inclusion : forall a b,
is_subcategory a b ->
ob_inj (subcategory_inclusion a b).

Lemma full_subcategory_inclusion : forall a b,
is_full_subcategory a b ->
full (subcategory_inclusion a b).

Lemma is_full_subcategory_rw : forall a b,
is_full_subcategory a b =
(is_subcategory a b & full (subcategory_inclusion a b)).

Definition ob_image_fs f :=
full_subcategory (target f) (fun x => inc x (ob_image f)).

Lemma is_subcategory_ob_image_fs : forall f,
Functor.axioms f ->
is_subcategory (ob_image_fs f) (target f).

Lemma is_full_subcategory_ob_image_fs : forall f,
Functor.axioms f ->
is_full_subcategory (ob_image_fs f) (target f).

Lemma inc_ob_image : forall f x,
Functor.axioms f ->
inc x (ob_image f) = (exists y, (ob (source f) y & fob f y = x)).

Lemma ob_ob_image_fs : forall f x,
Functor.axioms f ->
ob (ob_image_fs f) x = inc x (ob_image f).

Lemma mor_ob_image_fs : forall f u,
Functor.axioms f ->
mor (ob_image_fs f) u =
(mor (target f) u & inc (source u) (ob_image f) &
inc (target u) (ob_image f)).

Definition ob_image_factor f :=
Functor.create (source f) (ob_image_fs f)
(fmor f).

Lemma source_ob_image_factor : forall f,
source (ob_image_factor f) = source f.

Lemma target_ob_image_factor : forall f,
target (ob_image_factor f) = ob_image_fs f.

Lemma fmor_ob_image_factor : forall f u,
Functor.axioms f -> mor (source f) u ->
fmor (ob_image_factor f) u = fmor f u.

Lemma fob_ob_image_factor : forall f x,
Functor.axioms f -> ob (source f) x ->
fob (ob_image_factor f) x = fob f x.

Lemma ob_image_factor_axioms : forall f,
Functor.axioms f ->
Functor.axioms (ob_image_factor f).

Lemma ob_image_factorization : forall f,
Functor.axioms f ->
f = fcompose
(subcategory_inclusion (ob_image_fs f) (target f))
(ob_image_factor f).

Lemma ob_surj_ob_image_factor : forall f,
Functor.axioms f ->
ob_surj (ob_image_factor f).

Lemma ob_inj_ob_image_factor : forall f,
Functor.axioms f ->
ob_inj (ob_image_factor f) = ob_inj f.

Lemma ob_iso_ob_image_factor : forall f,
Functor.axioms f ->
ob_iso (ob_image_factor f) = ob_inj f.

Lemma full_ob_image_factor : forall f,
Functor.axioms f ->
full (ob_image_factor f) = full f.

Lemma faithful_ob_image_factor : forall f,
Functor.axioms f ->
faithful (ob_image_factor f) = faithful f.

Lemma iso_to_full_subcategory_rw : forall f,
iso_to_full_subcategory f =
(Functor.axioms f & has_finverse (ob_image_factor f)).

Lemma subcategory_extensionality : forall a b c,
is_subcategory a c -> is_subcategory b c ->
(forall u, mor a u -> mor b u) ->
(forall u, mor b u -> mor a u) ->
a = b.

Lemma full_subcategory_extensionality : forall a b c,
is_full_subcategory a c -> is_full_subcategory b c ->
(forall x, ob a x -> ob b x) ->
(forall x, ob b x -> ob a x) -> a = b.

Lemma iso_to_full_subcategory_interp :forall f,
iso_to_full_subcategory f =
(exists g, (is_full_subcategory (target g) (target f) &
has_finverse g &
f = fcompose (subcategory_inclusion (target g) (target f)) g)).

Definition mor_image f :=
Z (morphisms (target f)) (msays f).

Lemma sub_mor_image : forall f,
sub (mor_image f) (morphisms (target f)).

Lemma inc_mor_image : forall f u,
Functor.axioms f ->
inc u (mor_image f) = (exists v, (mor (source f) v &
fmor f v = u)).

Definition add_inverses a s :=
Z (morphisms a)
(fun y => inc y s \/ (invertible a y & inc (inverse a y) s)).

Lemma sub_add_inverses_refl : forall a s,
sub s (morphisms a) -> sub s (add_inverses a s).

Lemma sub_add_inverses_morphisms : forall a s,
sub (add_inverses a s) (morphisms a).

Lemma inc_add_inverses : forall a s y,
Category.axioms a -> sub s (morphisms a) ->
inc y (add_inverses a s) =
(mor a y & (inc y s \/ (invertible a y & inc (inverse a y) s))).

Definition inverse_closed a b :=
is_subcategory a b &
(forall u, mor a u -> invertible b u -> mor a (inverse b u)).

Lemma sub_add_inverses : forall a b s,
inverse_closed a b -> sub s (morphisms a) ->
sub (add_inverses b s) (morphisms a).

Definition generates a s :=
Category.axioms a & sub s (morphisms a) &
(forall b, is_subcategory b a -> sub s (morphisms b) -> b = a).

Lemma faithful_pull_morphism_criterion : forall f a,
Category.axioms a -> ob_surj f ->
faithful (pull_morphism a f).

Definition globular t :=
(Nat_Trans.like t) &
(Category.axioms (osource t)) &
(Category.axioms (otarget t)) &
(Functor.axioms (source t)) &
(Functor.axioms (target t)) &
(source (target t)) = (osource t) &
(target (source t)) = (otarget t).

Definition natural_ob t x :=
ob (osource t) x &
mor (otarget t) (ntrans t x) &
source (ntrans t x) = fob (source t) x &
target (ntrans t x) = fob (target t) x.

Definition natural_mor t u :=
mor (osource t) u &
natural_ob t (source u) &
natural_ob t (target u) &
comp (otarget t) (ntrans t (target u)) (fmor (source t) u)
= comp (otarget t) (fmor (target t) u) (ntrans t (source u)).

Lemma nat_trans_axioms_rw : forall t,
Nat_Trans.axioms t =
(globular t &
(forall x, ob (osource t) x -> natural_ob t x) &
(forall u, mor (osource t) u -> natural_mor t u)).

Lemma natural_mor_comp :forall t u v,
globular t -> natural_mor t u -> natural_mor t v ->
source u = target v -> natural_mor t (comp (osource t) u v).

Lemma natural_mor_id : forall t x,
globular t -> natural_ob t x ->
natural_mor t (id (osource t) x).

Lemma natural_ob_source :forall t u,
globular t -> natural_mor t u ->
natural_ob t (source u).

Lemma natural_ob_target : forall t u,
globular t -> natural_mor t u ->
natural_ob t (target u).

Definition naturality_subcategory t :=
subcategory (osource t) (natural_ob t) (natural_mor t).

Lemma is_subcategory_naturality_subcategory : forall t,
globular t -> is_subcategory (naturality_subcategory t) (osource t).

Lemma ob_naturality_subcategory : forall t x,
globular t -> ob (naturality_subcategory t) x = natural_ob t x.

Lemma mor_naturality_subcategory : forall t u,
globular t -> mor (naturality_subcategory t) u = natural_mor t u.

Lemma invertible_left_multiply : forall a u v w,
invertible a u -> mor a v -> mor a w ->
source u = target v -> source u = target w ->
comp a u v = comp a u w -> v = w.

Lemma natural_mor_inverse : forall t u,
globular t -> invertible (osource t) u ->
natural_mor t u -> natural_mor t (inverse (osource t) u).

Lemma inverse_closed_naturality_subcategory : forall t,
globular t -> inverse_closed (naturality_subcategory t) (osource t).

Lemma sub_add_inverses_naturality_subcategory : forall s t,
globular t -> sub s (morphisms (naturality_subcategory t)) ->
sub (add_inverses (osource t) s) 
(morphisms (naturality_subcategory t)).

Definition equalizer_subcategory f g :=
subcategory (source f) (fun x => (fob f x = fob g x)) 
(fun u => (fmor f u = fmor g u)).

Lemma equalizer_subcategory_property : forall f g,
Functor.axioms f -> Functor.axioms g ->
source f = source g -> target f = target g ->
subcategory_property (source f)
(fun x => (fob f x = fob g x)) (fun u => (fmor f u = fmor g u)).

Lemma is_subcategory_equalizer_subcategory : forall f g,
Functor.axioms f -> Functor.axioms g ->
source f = source g -> target f = target g ->
is_subcategory (equalizer_subcategory f g) (source f).

Lemma ob_equalizer_subcategory : forall f g x,
Functor.axioms f -> Functor.axioms g -> 
source f = source g -> target f = target g ->
ob (equalizer_subcategory f g) x = 
(ob (source f) x & fob f x = fob g x).

Lemma mor_equalizer_subcategory : forall f g u,
Functor.axioms f -> Functor.axioms g -> 
source f = source g -> target f = target g ->
mor (equalizer_subcategory f g) u = 
(mor (source f) u & fmor f u = fmor g u).

Lemma mor_equ_subcat_inverse : forall f g u,
Functor.axioms f -> Functor.axioms g -> 
source f = source g -> target f = target g ->
invertible (source f) u ->
mor (equalizer_subcategory f g) u ->
mor (equalizer_subcategory f g) (inverse (source f) u).

Lemma inverse_closed_equalizer_subcatgory : forall f g,
Functor.axioms f -> Functor.axioms g -> 
source f = source g -> target f = target g ->
inverse_closed (equalizer_subcategory f g) (source f).

Lemma mor_equ_subcat_fmor : forall f g h u,
Functor.axioms f -> Functor.axioms g -> Functor.axioms h ->
source g = target f -> source h = target f ->
fcompose g f = fcompose h f ->
mor (source f) u ->
mor (equalizer_subcategory g h) (fmor f u).

Lemma sub_mor_image_equalizer_subcategory : forall f g h,
Functor.axioms f -> Functor.axioms g -> Functor.axioms h ->
source g = target f -> source h = target f ->
fcompose g f = fcompose h f ->
sub (mor_image f) (morphisms (equalizer_subcategory g h)).

Lemma sub_add_inverses_equalizer_subcategory : forall f g h,
Functor.axioms f -> Functor.axioms g -> Functor.axioms h ->
source g = target f -> source h = target f ->
fcompose g f = fcompose h f ->
sub (add_inverses (target f) (mor_image f)) 
(morphisms (equalizer_subcategory g h)).
   
Lemma equalizer_subcategory_extensionality_criterion :
forall f g,
Functor.axioms f -> Functor.axioms g ->
source f = source g -> target f = target g ->
equalizer_subcategory f g = source f ->
f = g.

Lemma fcompose_eq_shows_eq : forall f g h,
Functor.axioms f -> Functor.axioms g -> Functor.axioms h ->
source g = target f -> source h = target f ->
fcompose g f = fcompose h f ->
generates (target f) (add_inverses (target f) (mor_image f)) ->
g = h.

Lemma ob_inj_pull_morphism_criterion : forall f a,
Category.axioms a ->
ob_iso f ->
generates (target f) (add_inverses (target f) (mor_image f)) ->
ob_inj (pull_morphism a f).

Definition full_pull_situation f g h u :=
Functor.axioms f &
Functor.axioms g &
Functor.axioms h &
Nat_Trans.axioms u &
source g = target f &
source h = target f &
fcompose g f = source u &
fcompose h f = target u &
ob_iso f &
generates (target f) (add_inverses (target f) (mor_image f)).

Lemma full_pull_additional_facts : forall f g h u,
full_pull_situation f g h u ->
(target g = otarget u &
target h = otarget u &
source f = osource u &
source g = source h &
target g = target h &
osource u = source f).

Definition ob_inverse_pr f x y:=
ob (source f) y & fob f y = x.

Definition ob_inverse f x :=
choose (ob_inverse_pr f x).

Lemma exists_ob_inverse_pr : forall f x,
ob_iso f -> ob (target f) x ->
exists y, ob_inverse_pr f x y.

Lemma fob_ob_inverse : forall f x,
ob_iso f -> ob (target f) x ->
fob f (ob_inverse f x) = x.

Lemma ob_ob_inverse : forall f x,
ob_iso f -> ob (target f) x ->
ob (source f) (ob_inverse f x).

Lemma ob_inverse_fob : forall f x,
ob_iso f -> ob (source f) x ->
ob_inverse f (fob f x) = x.

Definition ntdotted f g h u :=
Nat_Trans.create g h (fun x => ntrans u (ob_inverse f x)).

Lemma source_ntdotted : forall f g h u,
source (ntdotted f g h u) = g.

Lemma target_ntdotted : forall f g h u,
target (ntdotted f g h u) = h.

Lemma osource_ntdotted : forall f g h u,
osource (ntdotted f g h u) = source g.

Lemma otarget_ntdotted : forall f g h u,
otarget (ntdotted f g h u) = target h.

Lemma globular_ntdotted : forall f g h u,
full_pull_situation f g h u ->
globular (ntdotted f g h u).

Lemma ntrans_ntdotted : forall f g h u x,
full_pull_situation f g h u ->
ob (target f) x ->
ntrans (ntdotted f g h u) x =
ntrans u (ob_inverse f x).

Lemma ntrans_ntdotted_fob : forall f g h u x,
full_pull_situation f g h u ->
ob (source f) x ->
ntrans (ntdotted f g h u) (fob f x)
= ntrans u x.

Lemma natural_ob_ntdotted : forall f g h u x,
full_pull_situation f g h u ->
ob (target f) x ->
natural_ob (ntdotted f g h u) x.

Lemma natural_mor_ntdotted_fmor : forall f g h u y,
full_pull_situation f g h u ->
mor (source f) y ->
natural_mor (ntdotted f g h u) (fmor f y).

Lemma sub_mor_image_naturality_subcategory : forall f g h u,
full_pull_situation f g h u ->
sub (mor_image f) 
(morphisms (naturality_subcategory (ntdotted f g h u))).

Lemma sub_add_inverses_naturality_subcategory_ntdotted :
forall f g h u,
full_pull_situation f g h u ->
sub (add_inverses (target f) (mor_image f))
(morphisms (naturality_subcategory (ntdotted f g h u))).

Lemma naturality_subcategory_ntdotted_all : forall f g h u,
full_pull_situation f g h u ->
naturality_subcategory (ntdotted f g h u) = target f.

Lemma ntdotted_axioms : forall f g h u,
full_pull_situation f g h u ->
Nat_Trans.axioms (ntdotted f g h u).

Lemma htrans_right_ntdotted : forall f g h u,
full_pull_situation f g h u ->
htrans_right (ntdotted f g h u) f = u.

Lemma full_pull_morphism_criterion : forall f a,
Category.axioms a ->
ob_iso f ->
generates (target f) (add_inverses (target f) (mor_image f)) ->
full (pull_morphism a f).


Lemma iso_to_full_subcategory_pull_morphism_criterion : forall f a,
Category.axioms a ->
ob_iso f ->
generates (target f) (add_inverses (target f) (mor_image f)) ->
iso_to_full_subcategory (pull_morphism a f).

End Ob_Iso_Functor.

\end{verbatim}

\subsection{The module {\tt Associating\_Quotient}}

\begin{verbatim}
Module Associating_Quotient.
Export Quotient_Category.
Export Quotient_Functor.


Definition rqcat a :=
Category.Notations.like a &
(forall u, is_mor a u -> Arrow.like u) &
(forall x, is_ob a x -> is_mor a (id a x)) &
(forall x, is_ob a x -> source (id a x) = x) &
(forall x, is_ob a x -> target (id a x) = x) &
(forall u, is_mor a u -> is_ob a (source u)) &
(forall u, is_mor a u -> is_ob a (target u)) &
(forall u v, is_mor a u -> is_mor a v -> source u = target v ->
is_mor a (comp a u v)) &
(forall u v, is_mor a u -> is_mor a v -> source u = target v ->
source (comp a u v) = source v) &
(forall u v, is_mor a u -> is_mor a v -> source u = target v ->
target (comp a u v) = target u).

Lemma rq_is_ob_source : forall a u,
rqcat a -> is_mor a u -> is_ob a (source u).

Lemma rq_is_ob_target : forall a u,
rqcat a -> is_mor a u -> is_ob a (target u).

Lemma rq_is_mor_id : forall a x,
rqcat a -> is_ob a x -> is_mor a (id a x).

Lemma rq_source_id : forall a x,
rqcat a -> is_ob a x -> source (id a x) = x.

Lemma rq_target_id : forall a x,
rqcat a -> is_ob a x -> target (id a x) = x.

Lemma rq_is_mor_comp : forall a u v,
rqcat a -> is_mor a u -> is_mor a v -> source u = target v ->
is_mor a (comp a u v).

Lemma rq_source_comp : forall a u v,
rqcat a -> is_mor a u -> is_mor a v -> source u = target v ->
source (comp a u v) = source v.

Lemma rq_target_comp : forall a u v,
rqcat a -> is_mor a u -> is_mor a v -> source u = target v ->
target (comp a u v) = target u.

Definition left_id_ok a :=
(forall u, is_mor a u -> comp a (id a (target u)) u = u).

Definition right_id_ok a :=
(forall u, is_mor a u -> comp a u (id a (source u)) = u).

Definition assoc_ok a :=
(forall u v w, is_mor a u -> is_mor a v -> is_mor a w ->
source u = target v -> source v = target w ->
comp a (comp a u v) w = comp a u (comp a v w)).

Lemma cat_axioms_rw_rq : forall a,
Category.axioms a =
(rqcat a & left_id_ok a & right_id_ok a & assoc_ok a).

Definition rqcat_equiv_rel a r :=
rqcat a &
is_equivalence_relation r &
(forall x y, related r x y -> is_mor a x) &
(forall x y, related r x y -> is_mor a y) &
(forall x y, related r x y -> source x = source y) &
(forall x y, related r x y -> target x = target y) &
(forall u, is_mor a u -> related r u u) &
(forall x y u v, related r x y -> related r u v ->
source x = target u ->
related r (comp a x u) (comp a y v)) &
(forall u, is_mor a u ->
related r (comp a (id a (target u)) u) u) &
(forall u, is_mor a u ->
related r (comp a u (id a (source u))) u) &
(forall u v w, is_mor a u -> is_mor a v -> is_mor a w ->
source u = target v -> source v = target w ->
related r (comp a (comp a u v) w) (comp a u (comp a v w))).

Lemma rq_inc_arrow_arrow_class : forall a r u v,
rqcat_equiv_rel a r -> is_mor a v ->
inc v (arrow (arrow_class r u)) = related r u v.

Definition rq_quotient_arrow a r u :=
rqcat_equiv_rel a r &
(exists y, is_mor a y & u = arrow_class r y).

Lemma rq_quotient_arrow_arrow_class : forall a r u,
rqcat_equiv_rel a r -> is_mor a u ->
rq_quotient_arrow a r (arrow_class r u).

Lemma rq_inc_quotient_morphisms : forall a r u,
rqcat_equiv_rel a r ->
inc u (quotient_morphisms a r) = rq_quotient_arrow a r u.

Lemma rq_related_arrow_class_eq : forall a r u v,
rqcat_equiv_rel a r ->
is_mor a u ->
related r u v =
(arrow_class r u = arrow_class r v).

Lemma rq_inc_arrow_class_refl : forall a r u,
rqcat_equiv_rel a r -> is_mor a u -> inc u (arrow (arrow_class r u)).

Lemma rq_nonempty_arrow : forall a r u,
rq_quotient_arrow a r u ->
nonempty (arrow u).

Lemma rq_inc_arrow_rep_arrow : forall a r u,
rq_quotient_arrow a r u ->
inc (arrow_rep u) (arrow u).

Lemma rq_related_arrow_rep_arrow_class : forall a r u,
rqcat_equiv_rel a r -> is_mor a u ->
related r u (arrow_rep (arrow_class r u)).

Lemma rq_source_arrow_rep : forall a r u,
rq_quotient_arrow a r u -> source (arrow_rep u) = source u.

Lemma rq_target_arrow_rep : forall a r u,
rq_quotient_arrow a r u -> target (arrow_rep u) = target u.

Lemma rq_arrow_class_arrow_rep : forall a r u,
rq_quotient_arrow a r u ->
arrow_class r (arrow_rep u) = u.

Lemma rq_inc_arrow_facts : forall a r u y,
rq_quotient_arrow a r u -> inc y (arrow u) ->
(is_mor a y & source y = source u & target y = target u).

Lemma rq_mor_arrow_rep : forall a r u,
rq_quotient_arrow a r u ->
is_mor a (arrow_rep u).

Lemma rq_related_arrow_rep : forall a r u v,
rqcat_equiv_rel a r -> is_mor a v -> u = arrow_class r v ->
related r v (arrow_rep u).

Lemma rq_related_arrow_rep_rw : forall a r u v,
rqcat_equiv_rel a r -> rq_quotient_arrow a r u -> is_mor a v ->
related r v (arrow_rep u) = (u = arrow_class r v).

Lemma rq_source_quot_comp : forall a r u v,
rqcat_equiv_rel a r ->
rq_quotient_arrow a r u ->
rq_quotient_arrow a r v -> source u = target v ->
source (quot_comp a r u v) = source v.

Lemma rq_target_quot_comp : forall a r u v,
rqcat_equiv_rel a r ->
rq_quotient_arrow a r u ->
rq_quotient_arrow a r v -> source u = target v ->
target (quot_comp a r u v) = target u.

Lemma rq_quotient_arrow_quot_id : forall a r x,
rqcat_equiv_rel a r -> is_ob a x ->
rq_quotient_arrow a r (quot_id a r x).

Lemma rq_quotient_arrow_quot_comp : forall a r u v,
rqcat_equiv_rel a r -> rq_quotient_arrow a r u ->
rq_quotient_arrow a r v
-> source u = target v ->
rq_quotient_arrow a r (quot_comp a r u v).

Lemma rq_arrow_class_eq : forall a r u v,
rqcat_equiv_rel a r -> is_mor a u ->
(arrow_class r u = arrow_class r v) =
(related r u v).

 
Lemma rq_related_comp : forall a r x y u v,
rqcat_equiv_rel a r ->
related r x y -> related r u v -> source x = target u ->
related r (comp a x u) (comp a y v).

Lemma rq_quot_comp_arrow_class : forall a r u v,
rqcat_equiv_rel a r -> is_mor a u -> is_mor a v ->
source u = target v ->
quot_comp a r (arrow_class r u) (arrow_class r v) =
arrow_class r (comp a u v).

Lemma rq_left_id_related : forall a r u,
rqcat_equiv_rel a r ->
is_mor a u ->
related r (comp a (id a (target u )) u) u.

Lemma rq_right_id_related : forall a r u,
rqcat_equiv_rel a r ->
is_mor a u ->
related r (comp a u (id a (source u ))) u.

Lemma rq_assoc_related : forall a r u v w,
rqcat_equiv_rel a r ->
is_mor a u -> is_mor a v -> is_mor a w ->
source u = target v -> source v = target w ->
related r (comp a (comp a u v) w) (comp a u (comp a v w)).

Lemma rq_quot_id_left : forall a r a' r' x u,
rqcat_equiv_rel a r -> rq_quotient_arrow a r u ->
x = target u -> a' = a -> r' = r ->
quot_comp a r (quot_id a' r' x) u = u.

Lemma rq_quot_id_right : forall a r a' r' x u,
rqcat_equiv_rel a r -> rq_quotient_arrow a r u ->
x = source u -> a' = a -> r' = r ->
quot_comp a r u (quot_id a' r' x) = u.

Lemma rq_quot_comp_assoc : forall a r a' r' u v w,
rqcat_equiv_rel a r -> rq_quotient_arrow a r u ->
rq_quotient_arrow a r v
-> rq_quotient_arrow a r w -> source u = target v ->
source v = target w -> a' = a -> r' = r ->
quot_comp a r (quot_comp a' r' u v) w = 
quot_comp a r u (quot_comp a' r' v w).

Lemma rq_is_mor_quotient_cat : forall a r u,
rqcat_equiv_rel a r ->
is_mor (quotient_cat a r) u = rq_quotient_arrow a r u.

Lemma rq_comp_quotient_cat : forall a r u v,
rqcat_equiv_rel a r -> rq_quotient_arrow a r u ->
rq_quotient_arrow a r v ->
source u = target v ->
comp (quotient_cat a r) u v = quot_comp a r u v.

Lemma rq_id_quotient_cat : forall a r x,
rqcat_equiv_rel a r -> is_ob a x ->
id (quotient_cat a r) x = quot_id a r x.

Lemma rq_source_quot_id : forall a r x,
rqcat_equiv_rel a r -> is_ob a x ->
source (quot_id a r x) = x.

Lemma rq_target_quot_id : forall a r x,
rqcat_equiv_rel a r -> is_ob a x ->
target (quot_id a r x) = x.

Lemma rqcat_quotient_cat : forall a r,
rqcat_equiv_rel a r ->
rqcat (quotient_cat a r).

Lemma rq_quotient_cat_axioms : forall a r,
rqcat_equiv_rel a r -> Category.axioms (quotient_cat a r).

Lemma rq_ob_quotient_cat : forall a r x,
rqcat_equiv_rel a r ->
ob (quotient_cat a r) x = is_ob a x.

Lemma rq_mor_quotient_cat : forall a r u,
rqcat_equiv_rel a r ->
mor (quotient_cat a r) u = rq_quotient_arrow a r u.

Lemma rq_mor_quotient_cat_ex : forall a r u,
rqcat_equiv_rel a r ->
mor (quotient_cat a r) u =
(exists y, is_mor a y & u = arrow_class r y).

Lemma rq_mor_quotient_cat_quot_id : forall a r x,
rqcat_equiv_rel a r -> is_ob a x ->
mor (quotient_cat a r) (quot_id a r x).

Lemma rq_mor_quotient_cat_quot_comp : forall a r u v,
rqcat_equiv_rel a r ->
mor (quotient_cat a r) u ->
mor (quotient_cat a r) v ->
source u = target v ->
mor (quotient_cat a r) (quot_comp a r u v).

Lemma rq_mor_quotient_cat_arrow_class : forall a r u,
rqcat_equiv_rel a r -> is_mor a u ->
mor (quotient_cat a r) (arrow_class r u).


Definition rqfunctor_property a b r fo fm :=
Category.axioms a &
rqcat_equiv_rel b r &
(forall x, ob a x -> is_ob b (fo x)) &
(forall u, mor a u -> is_mor b (fm u)) &
(forall u, mor a u -> source (fm u) = fo (source u)) &
(forall u, mor a u -> target (fm u) = fo (target u)) &
(forall x, ob a x -> related r (fm (id a x)) (id b (fo x))) &
(forall u v, mor a u -> mor a v -> source u = target v ->
related r (fm (comp a u v)) (comp b (fm u) (fm v))).

Lemma rq_fob_qfunctor : forall a b r fo fm x,
rqfunctor_property a b r fo fm -> ob a x ->
fob (qfunctor a b r fo fm) x = fo x.

Lemma rq_ob_fob_qfunctor : forall a b r fo fm x,
rqfunctor_property a b r fo fm -> ob a x ->
ob (quotient_cat b r) (fob (qfunctor a b r fo fm) x).

Lemma rq_mor_fmor_qfunctor : forall a b r fo fm u,
rqfunctor_property a b r fo fm -> mor a u ->
mor (quotient_cat b r) (fmor (qfunctor a b r fo fm) u).

Lemma rq_source_fmor_qfunctor : forall a b r fo fm u,
rqfunctor_property a b r fo fm -> mor a u ->
source (fmor (qfunctor a b r fo fm) u) = 
fob (qfunctor a b r fo fm) (source u).

Lemma rq_target_fmor_qfunctor : forall a b r fo fm u,
rqfunctor_property a b r fo fm -> mor a u ->
target (fmor (qfunctor a b r fo fm) u) = 
fob (qfunctor a b r fo fm) (target u).

Lemma rq_fmor_qfunctor_id : forall a b r fo fm x,
rqfunctor_property a b r fo fm -> ob a x ->
fmor (qfunctor a b r fo fm) (id a x) =
id (quotient_cat b r) (fob (qfunctor a b r fo fm) x).

Lemma rq_fmor_qfunctor_comp : forall a b r fo fm u v,
rqfunctor_property a b r fo fm -> mor a u ->
mor a v -> source u = target v ->
fmor (qfunctor a b r fo fm) (comp a u v) =
comp (quotient_cat b r) (fmor (qfunctor a b r fo fm) u)
(fmor (qfunctor a b r fo fm) v).

Lemma qfunctor_axioms : forall a b r fo fm,
rqfunctor_property a b r fo fm ->
Functor.axioms (qfunctor a b r fo fm).

End Associating_Quotient.

\end{verbatim}

\newpage

\section{The file {\tt gzdef.v}}

\begin{verbatim}

Require Export freecat.
Require Export qcat.

\end{verbatim}

\subsection{The module {\tt GZ\_Def}}

\begin{verbatim}
Module GZ_Def.
Export Free_Category.
Export Quotient_Functor.

Definition Forward := R (f_(o_(r_ DOT))).
Definition Backward := R (b_(k_(d_ DOT))).

Definition forward_arrow u :=
Arrow.create (source u) (target u) (pair Forward u).

Definition backward_arrow u :=
Arrow.create (target u) (source u) (pair Backward u).

Lemma source_forward_arrow : forall u,
source (forward_arrow u) = source u.

Lemma target_forward_arrow : forall u,
target (forward_arrow u) = target u.

Lemma source_backward_arrow : forall u,
source (backward_arrow u) = target u.

Lemma target_backward_arrow : forall u,
target (backward_arrow u) = source u.

Definition original_arrow u := pr2 (arrow u).

Lemma original_arrow_forward_arrow : forall u,
original_arrow (forward_arrow u) = u.

Lemma original_arrow_backward_arrow : forall u,
original_arrow (backward_arrow u) = u.

Definition direction u := pr1 (arrow u).

Lemma direction_forward_arrow : forall u,
direction (forward_arrow u) = Forward.

Lemma direction_backward_arrow : forall u,
direction (backward_arrow u) = Backward.

Definition localizing_system a s :=
Category.axioms a &
(forall u, inc u s -> mor a u).

Definition loc_edges a s :=
union2 (Image.create (morphisms a) forward_arrow) 
(Image.create s backward_arrow).

Lemma inc_loc_edges : forall a s u,
localizing_system a s ->
inc u (loc_edges a s) =
((exists y, (mor a y & u = forward_arrow y)) \/
(exists y, (mor a y & inc y s & u = backward_arrow y))).

Definition gz_graph a s := Graph.create (objects a) (loc_edges a s).

Lemma inc_vertices_gz_graph : forall a s x,
localizing_system a s ->
inc x (vertices (gz_graph a s)) = ob a x.

Lemma inc_edges_gz_graph : forall a s u,
localizing_system a s ->
inc u (edges (gz_graph a s)) = inc u (loc_edges a s).

Lemma gz_graph_axioms : forall a s,
localizing_system a s -> Graph.axioms (gz_graph a s).

Definition gz_freecat a s := freecat (gz_graph a s).

Lemma gz_freecat_axioms : forall a s,
localizing_system a s -> Category.axioms (gz_freecat a s).

Lemma ob_gz_freecat : forall a s x,
localizing_system a s ->
ob (gz_freecat a s) x = ob a x.

Lemma mor_gz_freecat : forall a s u,
localizing_system a s ->
mor (gz_freecat a s) u = mor_freecat (gz_graph a s) u.

Lemma comp_gz_freecat : forall a s u v,
localizing_system a s ->
mor (gz_freecat a s) u ->
mor (gz_freecat a s) v ->
source u = target v ->
comp (gz_freecat a s) u v = freecat_comp u v.

Lemma id_gz_freecat : forall a s x,
localizing_system a s ->
ob a x -> id (gz_freecat a s) x = freecat_id x.

Definition forward_edge u :=
freecat_edge (forward_arrow u).

Definition backward_edge u :=
freecat_edge (backward_arrow u).

Lemma source_forward_edge : forall u,
source (forward_edge u) = source u.

Lemma target_forward_edge : forall u,
target (forward_edge u) = target u.

Lemma source_backward_edge : forall u,
source (backward_edge u) = target u.

Lemma target_backward_edge : forall u,
target (backward_edge u) = source u.

Lemma inc_forward_arrow_loc_edges : forall a s u,
localizing_system a s ->mor a u ->
inc (forward_arrow u) (loc_edges a s).

Lemma inc_backward_arrow_loc_edges : forall a s q,
localizing_system a s ->inc q s ->
inc (backward_arrow q) (loc_edges a s).

Lemma mor_forward_edge : forall a s u,
localizing_system a s ->
mor a u -> mor (gz_freecat a s) (forward_edge u).

Lemma mor_backward_edge : forall a s q,
localizing_system a s ->
inc q s -> mor (gz_freecat a s) (backward_edge q).

Definition gz_rel a s :=
Z (coarse (gz_freecat a s))
(fun z =>
((exists x, (ob a x &
z = pair (forward_edge (id a x)) (freecat_id x))) \/
(exists q, (inc q s & z = pair
(freecat_comp (forward_edge q) (backward_edge q))
(freecat_id (target q)))) \/
(exists q, (inc q s & z = pair
(freecat_comp (backward_edge q) (forward_edge q))
(freecat_id (source q)))) \/
(exists u, exists v, (mor a u & mor a v & source u = target v &
z = pair (freecat_comp (forward_edge u) (forward_edge v))
(forward_edge (comp a u v)))))).

Lemma inc_coarse : forall a z,
Category.axioms a ->
inc z (coarse a) =
(exists u, exists v, (mor a u & mor a v & source u = source v
& target u = target v & z=pair u v)).

Lemma inc_gz_rel : forall a s z,
localizing_system a s ->
inc z (gz_rel a s) =
((exists x, (ob a x &
z = pair (forward_edge (id a x)) (freecat_id x)))\/
(exists q, (inc q s & z = pair
(freecat_comp (forward_edge q) (backward_edge q))
(freecat_id (target q))))\/
(exists q, (inc q s & z = pair
(freecat_comp (backward_edge q) (forward_edge q))
(freecat_id (source q))))\/
(exists u, exists v, (mor a u & mor a v & source u = target v &
z = pair (freecat_comp (forward_edge u) (forward_edge v))
(forward_edge (comp a u v))))).

Lemma related_gz_rel : forall a s e f,
localizing_system a s ->
related (gz_rel a s) e f =
((exists x, (ob a x &
e = (forward_edge (id a x)) & f = (freecat_id x))) \/
(exists q, (inc q s & e =
(freecat_comp (forward_edge q) (backward_edge q))
& f = (freecat_id (target q)))) \/
(exists q, (inc q s & e =
(freecat_comp (backward_edge q) (forward_edge q))
& f = (freecat_id (source q)))) \/
(exists u, exists v, (mor a u & mor a v & source u = target v &
e = (freecat_comp (forward_edge u) (forward_edge v))
& f = (forward_edge (comp a u v))))).

Lemma sub_gz_rel_coarse : forall a s,
sub (gz_rel a s) (coarse (gz_freecat a s)).

Lemma cat_rel_gz_rel : forall a s,
localizing_system a s ->
cat_rel (gz_freecat a s) (gz_rel a s).

Definition gz_cer a s := cer (gz_freecat a s) (gz_rel a s).

Lemma cat_equiv_rel_gz_cer : forall a s,
localizing_system a s ->
cat_equiv_rel (gz_freecat a s) (gz_cer a s).

Lemma related_gz_cer_first_mor : forall a s e f,
localizing_system a s ->
related (gz_cer a s) e f -> mor (gz_freecat a s) e.

Lemma related_gz_cer_second_mor : forall a s e f,
localizing_system a s ->
related (gz_cer a s) e f -> mor (gz_freecat a s) f.

Lemma related_gz_cer_same_source : forall a s e f,
localizing_system a s ->
related (gz_cer a s) e f -> source e = source f.

Lemma related_gz_cer_same_target : forall a s e f,
localizing_system a s ->
related (gz_cer a s) e f -> target e = target f.

Lemma related_cer : forall a r u v,
cat_rel a r -> related r u v -> related (cer a r) u v.

Lemma related_gz_cer_forward_edge_id :
forall a s x,
localizing_system a s ->
ob a x -> related (gz_cer a s) (forward_edge (id a x)) (freecat_id x).

Lemma related_gz_cer_comp_forward_backward : forall a s q,
localizing_system a s -> inc q s ->
related (gz_cer a s) (freecat_comp (forward_edge q) (backward_edge q))
(freecat_id (target q)).

Lemma related_gz_cer_comp_backward_forward : forall a s q,
localizing_system a s -> inc q s ->
related (gz_cer a s) (freecat_comp (backward_edge q) (forward_edge q))
(freecat_id (source q)).

Lemma related_gz_cer_comp_forward : forall a s u v,
localizing_system a s -> mor a u -> mor a v ->
source u = target v ->
related (gz_cer a s) (freecat_comp (forward_edge u) (forward_edge v))
(forward_edge (comp a u v)).

Lemma compatible_gz_cer_criterion : forall a s f,
localizing_system a s ->
Functor.axioms f ->
source f = gz_freecat a s ->
(forall x, ob a x -> fmor f (forward_edge (id a x)) = 
id (target f) (fob f x)) ->
(forall q, inc q s ->
are_inverse (target f) (fmor f (forward_edge q)) 
(fmor f (backward_edge q))) ->
(forall u v, mor a u -> mor a v -> source u = target v ->
comp (target f) (fmor f (forward_edge u)) (fmor f (forward_edge v)) =
fmor f (forward_edge (comp a u v))) ->
compatible (gz_cer a s) f.

Definition gz_loc a s :=
quotient_cat (gz_freecat a s) (gz_cer a s).

Lemma gz_loc_axioms : forall a s,
localizing_system a s ->
Category.axioms (gz_loc a s).

Definition gz_qprojection a s :=
qprojection (gz_freecat a s) (gz_cer a s).

Lemma source_gz_qprojection : forall a s,
localizing_system a s ->
source (gz_qprojection a s) = gz_freecat a s.

Lemma target_gz_qprojection : forall a s,
localizing_system a s ->
target (gz_qprojection a s) = gz_loc a s.

Lemma gz_qprojection_axioms : forall a s,
localizing_system a s ->
Functor.axioms (gz_qprojection a s).

Definition gz_proj a s :=
qfunctor a (gz_freecat a s) (gz_cer a s)
(fun x => x) forward_edge.

Lemma source_gz_proj : forall a s,
localizing_system a s ->
source (gz_proj a s) = a.

Lemma target_gz_proj : forall a s,
localizing_system a s ->
target (gz_proj a s) = gz_loc a s.

Definition gz_forward a s u :=
arrow_class (gz_cer a s) (forward_edge u).

Definition gz_backward a s u :=
arrow_class (gz_cer a s) (backward_edge u).

Lemma source_gz_forward : forall a s u,
localizing_system a s ->
mor a u -> source (gz_forward a s u) = source u.

Lemma target_gz_forward : forall a s u,
localizing_system a s ->
mor a u -> target (gz_forward a s u) = target u.

Lemma source_gz_backward : forall a s q,
localizing_system a s ->
inc q s -> source (gz_backward a s q) = target q.

Lemma target_gz_backward : forall a s q,
localizing_system a s ->
inc q s -> target (gz_backward a s q) = source q.

Lemma mor_gz_forward : forall a s u,
localizing_system a s ->
mor a u -> mor (gz_loc a s) (gz_forward a s u).

Lemma mor_gz_backward : forall a s q,
localizing_system a s ->
inc q s -> mor (gz_loc a s) (gz_backward a s q).

Lemma qfunctor_property_forward_edge :
forall a s, localizing_system a s ->
qfunctor_property a (gz_freecat a s) (gz_cer a s)
(fun x => x) forward_edge.

Lemma gz_proj_axioms : forall a s,
localizing_system a s ->
Functor.axioms (gz_proj a s).

Lemma fob_gz_proj : forall a s x,
localizing_system a s ->
ob a x ->
fob (gz_proj a s) x = x.

Lemma fmor_gz_proj : forall a s u,
localizing_system a s ->
mor a u -> fmor (gz_proj a s) u =
gz_forward a s u.

Lemma gz_forward_id : forall a s x,
localizing_system a s ->
ob a x -> gz_forward a s (id a x) = id (gz_loc a s) x.

Lemma comp_gz_forward : forall a s u v,
localizing_system a s -> mor a u -> mor a v ->
source u = target v ->
comp (gz_loc a s) (gz_forward a s u) (gz_forward a s v) =
gz_forward a s (comp a u v).

Lemma are_inverse_gz_forward_backward : forall a s q,
localizing_system a s ->
inc q s -> are_inverse (gz_loc a s) (gz_forward a s q)
(gz_backward a s q).

Lemma invertible_gz_loc_gz_forward : forall a s q,
localizing_system a s -> inc q s ->
invertible (gz_loc a s) (gz_forward a s q).

End GZ_Def.

\end{verbatim}

\subsection{The module {\tt GZ\_Thm}}

\begin{verbatim}
Module GZ_Thm.
Export Ob_Iso_Functor.
Export GZ_Def.


Lemma inc_gz_forward_mor_image : forall a s u,
localizing_system a s ->
mor a u -> inc (gz_forward a s u) (mor_image (gz_proj a s)).

Lemma inc_gz_backward_add_inverses :
forall a s q,
localizing_system a s ->
inc q s -> inc (gz_backward a s q) 
(add_inverses (gz_loc a s) (mor_image (gz_proj a s))).

Lemma gz_loc_induction : forall a s (P:E->Prop),
localizing_system a s ->
(forall u, mor a u -> P (gz_forward a s u)) ->
(forall q, inc q s -> P (gz_backward a s q)) ->
(forall x, ob a x -> P (id (gz_loc a s) x)) ->
(forall u v, mor (gz_loc a s) u ->
mor (gz_loc a s) v -> source u = target v ->
P u -> P v -> P (comp (gz_loc a s) u v)) ->
(forall y, mor (gz_loc a s) y -> P y).

Lemma gz_loc_subcategory_all_criterion : forall a s b,
localizing_system a s -> is_subcategory b (gz_loc a s) ->
(forall u, mor a u -> mor b (gz_forward a s u)) ->
(forall q, inc q s -> mor b (gz_backward a s q)) ->
b = gz_loc a s.

Lemma add_inverses_mor_image_gz_proj_generates_gz_loc : forall a s,
localizing_system a s ->
generates (gz_loc a s) 
(add_inverses (gz_loc a s) (mor_image (gz_proj a s))).

Lemma ob_gz_loc : forall a s x,
localizing_system a s ->
ob (gz_loc a s) x = ob a x.

Lemma ob_iso_gz_proj : forall a s,
localizing_system a s -> ob_iso (gz_proj a s).


Lemma iso_to_subcategory_pull_gz_proj : forall a s b,
localizing_system a s -> Category.axioms b ->
iso_to_full_subcategory (pull_morphism b (gz_proj a s)).


Lemma gz_proj_epimorphic : forall a s f g,
localizing_system a s -> Functor.axioms f ->
Functor.axioms g -> source f = gz_loc a s ->
source g = gz_loc a s ->
fcompose f (gz_proj a s) = fcompose g (gz_proj a s) ->
f = g.


Definition loc_compatible a s f :=
localizing_system a s &
Functor.axioms f &
source f = a &
(forall q, inc q s -> invertible (target f) (fmor f q)).

Definition fr_dotted a s f :=
free_functor (gz_graph a s) (target f)
(fun x => fob f x)
(fun u => Y (direction u = Forward) (fmor f (original_arrow u))
(inverse (target f) (fmor f (original_arrow u)))).

Lemma backward_neq_forward : Backward <> Forward.

Lemma Y_etc_forward_arrow : forall f u,
Y (direction (forward_arrow u) = Forward) 
(fmor f (original_arrow (forward_arrow u)))
(inverse (target f)
(fmor f (original_arrow (forward_arrow u)))) = fmor f u.

Lemma Y_etc_backward_arrow : forall f q,
Y (direction (backward_arrow q) = Forward)
(fmor f (original_arrow (backward_arrow q)))
(inverse (target f)
(fmor f (original_arrow (backward_arrow q)))) =
inverse (target f) (fmor f q).

Lemma fr_dotted_property : forall a s f,
loc_compatible a s f -> free_functor_property
(gz_graph a s) (target f)
(fun x => fob f x)
(fun u => Y (direction u = Forward) (fmor f (original_arrow u))
(inverse (target f) (fmor f (original_arrow u)))).

Lemma source_fr_dotted : forall a s f,
source (fr_dotted a s f) = gz_freecat a s.

Lemma target_fr_dotted : forall a s f,
target (fr_dotted a s f) = target f.

Lemma fob_fr_dotted : forall a s f x,
loc_compatible a s f -> ob a x ->
fob (fr_dotted a s f) x = fob f x.

Lemma fmor_fr_dotted_forward_edge : forall a s f u,
loc_compatible a s f -> mor a u ->
fmor (fr_dotted a s f) (forward_edge u) =
fmor f u.

Lemma fmor_fr_dotted_backward_edge : forall a s f q,
loc_compatible a s f -> inc q s ->
fmor (fr_dotted a s f) (backward_edge q) =
inverse (target f) (fmor f q).

Lemma fr_dotted_axioms : forall a s f,
loc_compatible a s f ->
Functor.axioms (fr_dotted a s f).

Lemma compatible_gz_cer_fr_dotted : forall a s f,
loc_compatible a s f ->
compatible (gz_cer a s) (fr_dotted a s f).

Definition gz_dotted a s f :=
qdotted (gz_cer a s) (fr_dotted a s f).

Lemma source_gz_dotted : forall a s f,
source (gz_dotted a s f) = gz_loc a s.

Lemma target_gz_dotted : forall a s f,
target (gz_dotted a s f) = target f.

Lemma fob_gz_dotted : forall a s f x,
loc_compatible a s f -> ob a x ->
fob (gz_dotted a s f) x = fob f x.

Lemma fmor_gz_dotted_gz_forward : forall a s f u,
loc_compatible a s f -> mor a u ->
fmor (gz_dotted a s f) (gz_forward a s u) = fmor f u.

Lemma fmor_gz_dotted_gz_backward : forall a s f q,
loc_compatible a s f -> inc q s ->
fmor (gz_dotted a s f) (gz_backward a s q) = 
inverse (target f) (fmor f q).

Lemma gz_dotted_axioms : forall a s f,
loc_compatible a s f ->
Functor.axioms (gz_dotted a s f).


Lemma fcompose_gz_dotted_gz_proj : forall a s f,
loc_compatible a s f ->
fcompose (gz_dotted a s f) (gz_proj a s) = f.


Lemma loc_compatible_fcompose : forall a s g,
localizing_system a s ->
Functor.axioms g ->
source g = gz_loc a s ->
loc_compatible a s (fcompose g (gz_proj a s)).


Lemma ob_image_pull_gz_proj : forall a s b f,
localizing_system a s ->
Category.axioms b ->
inc f (ob_image (pull_morphism b (gz_proj a s))) =
(Functor.axioms f & source f = a & target f = b &
(forall q, inc q s -> invertible b (fmor f q))).

End GZ_Thm.

\end{verbatim}

\newpage

\section{The file {\tt left\_fractions.v}}

\begin{verbatim}

Require Export gzdef.

\end{verbatim}

\subsection{The module {\tt Left\_Fractions}}

\begin{verbatim}
Module Left_Fractions.
Export GZ_Def.

Definition lf_symbol f t :=
Arrow.create (source f) (source t) (pair f t).

Definition lf_forward v := pr1 (arrow v).
Definition lf_backward v := pr2 (arrow v).

Lemma source_lf_symbol : forall f t, source (lf_symbol f t)
= source f.

Lemma target_lf_symbol : forall f t, target (lf_symbol f t)
= source t.

Lemma lf_forward_lf_symbol : forall f t, lf_forward
(lf_symbol f t) = f.

Lemma lf_backward_lf_symbol : forall f t, lf_backward
(lf_symbol f t) = t.

Definition lf_symbol_like v :=
v = lf_symbol (lf_forward v) (lf_backward v).

Lemma lf_symbol_like_lf_symbol : forall f t,
lf_symbol_like (lf_symbol f t).

Definition lf_symbol_property a s f t :=
localizing_system a s &
mor a f & inc t s & target f = target t.

Definition is_lf_symbol a s v :=
lf_symbol_like v & lf_symbol_property a s (lf_forward v)
(lf_backward v).

Lemma is_lf_symbol_lf_symbol : forall a s f t,
lf_symbol_property a s f t ->
is_lf_symbol a s (lf_symbol f t).

Definition multiplicative_system a s :=
localizing_system a s &
(forall y z, inc y s -> inc z s ->
source y = target z -> inc (comp a y z) s).

Definition has_left_fractions a s :=
multiplicative_system a s
&
(forall x, ob a x -> inc (id a x) s)
&
(forall r g, inc r s -> mor a g -> source r = source g ->
exists u, (is_lf_symbol a s u & source u = target r
& target u = target g & comp a (lf_backward u) g =
comp a (lf_forward u) r))
&
(forall v r t, inc v s -> mor a r -> mor a t ->
source r = target v -> source t = target v ->
comp a r v = comp a t v ->
exists w, (inc w s & source w = target r & source w = target t
& comp a w r = comp a w t)).

Definition lf_choice a s r g :=
choose (fun u => (is_lf_symbol a s u & source u = target r
& target u = target g & comp a (lf_backward u) g =
comp a (lf_forward u) r)).

Lemma is_lf_symbol_lf_choice : forall a s r g,
has_left_fractions a s ->
inc r s -> mor a g -> source r = source g ->
is_lf_symbol a s (lf_choice a s r g).

Lemma source_lf_choice : forall a s r g,
has_left_fractions a s ->
inc r s -> mor a g -> source r = source g ->
source (lf_choice a s r g) = target r.

Lemma target_lf_choice : forall a s r g,
has_left_fractions a s ->
inc r s -> mor a g -> source r = source g ->
target (lf_choice a s r g) = target g.

Lemma comp_lf_backward_lf_choice : forall a s r g,
has_left_fractions a s ->
inc r s -> mor a g -> source r = source g ->
comp a (lf_backward (lf_choice a s r g)) g
 = comp a (lf_forward (lf_choice a s r g)) r.

Lemma source_lf_forward : forall v,
lf_symbol_like v -> source (lf_forward v) = source v.

Lemma source_lf_backward : forall v,
lf_symbol_like v -> source (lf_backward v) = target v.

Lemma inc_ms_mor : forall a u,
(exists s, (localizing_system a s & inc u s)) ->
mor a u.

Lemma inc_hlf_mor : forall a u,
(exists s, (has_left_fractions a s & inc u s)) ->
mor a u.

Definition lf_id_rep a x :=
lf_symbol (id a x) (id a x).

Lemma source_lf_id_rep : forall a x,
ob a x -> source (lf_id_rep a x) = x.

Lemma target_lf_id_rep : forall a x,
ob a x -> target (lf_id_rep a x) = x.

Lemma is_lf_symbol_lf_id_rep : forall a s x,
has_left_fractions a s -> ob a x ->
is_lf_symbol a s (lf_id_rep a x).

Definition lf_vertex v :=
target (lf_backward v).

Lemma target_lf_forward : forall u,
(exists a, exists s, (is_lf_symbol a s u)) ->
target (lf_forward u) = lf_vertex u.

Lemma target_lf_backward : forall u,
target (lf_backward u) = lf_vertex u.

Definition lf_extend a (s:E) p u :=
lf_symbol (comp a p (lf_forward u)) (comp a p (lf_backward u)).

Lemma mor_lf_forward : forall a u,
(exists s, is_lf_symbol a s u) ->
mor a (lf_forward u).

Lemma inc_lf_backward : forall s u,
(exists a, is_lf_symbol a s u) ->
inc (lf_backward u) s.

Lemma mor_lf_backward : forall a u,
(exists s, is_lf_symbol a s u) ->
mor a (lf_backward u).

Lemma is_lf_symbol_lf_extend : forall a s p u,
has_left_fractions a s ->
is_lf_symbol a s u ->
mor a p -> source p = lf_vertex u ->
inc (comp a p (lf_backward u)) s ->
is_lf_symbol a s (lf_extend a s p u).

Lemma source_lf_extend : forall a s p u,
is_lf_symbol a s u ->
mor a p -> source p = lf_vertex u ->
source (lf_extend a s p u) = source u.

Lemma target_lf_extend : forall a s p u,
is_lf_symbol a s u ->
mor a p -> source p = lf_vertex u ->
target (lf_extend a s p u) = target u.

Lemma lf_forward_lf_extend : forall a s p u,
lf_forward (lf_extend a s p u) = comp a p (lf_forward u).

Lemma lf_backward_lf_extend : forall a s p u,
lf_backward (lf_extend a s p u) = comp a p (lf_backward u).

Lemma lf_vertex_lf_extend : forall a s p u,
has_left_fractions a s ->
is_lf_symbol a s u -> mor a p -> source p = lf_vertex u ->
lf_vertex (lf_extend a s p u) = target p.


Definition lf_beyond a s u v :=
has_left_fractions a s &
is_lf_symbol a s u & is_lf_symbol a s v &
source u = source v & target u = target v &
(exists p, (mor a p & source p = lf_vertex u &
lf_extend a s p u = v)).

Definition lf_under a s u v :=
has_left_fractions a s &
is_lf_symbol a s u & is_lf_symbol a s v &
source u = source v & target u = target v &
(exists p, (inc p s & source p = lf_vertex u &
lf_extend a s p u = v)).

Lemma lf_under_lf_beyond : forall a s u v,
lf_under a s u v -> lf_beyond a s u v.

Lemma lf_beyond_lf_extend : forall a s p u,
has_left_fractions a s ->
is_lf_symbol a s u -> mor a p ->
source p = lf_vertex u ->
inc (comp a p (lf_backward u)) s ->
lf_beyond a s u (lf_extend a s p u).

Lemma inc_comp : forall a s u v,
multiplicative_system a s ->
inc u s -> inc v s -> source u = target v ->
inc (comp a u v) s.

Lemma lf_under_lf_extend : forall a s p u,
has_left_fractions a s ->
is_lf_symbol a s u -> inc p s ->
source p = lf_vertex u ->
lf_under a s u (lf_extend a s p u).

Lemma lf_symbol_like_extens : forall u v,
lf_symbol_like u -> lf_symbol_like v ->
lf_forward u = lf_forward v -> lf_backward u = lf_backward v ->
u = v.

Lemma lf_symbol_like_lf_extend : forall a s p u,
lf_symbol_like (lf_extend a s p u).

Lemma lf_extend_comp : forall a s u p q,
is_lf_symbol a s u ->
mor a p -> mor a q ->
source q = lf_vertex u ->
source p = target q ->
lf_extend a s (comp a p q) u =
lf_extend a s p (lf_extend a s q u).

Lemma lf_beyond_trans : forall a s u v w,
lf_beyond a s u v -> lf_beyond a s v w ->
lf_beyond a s u w.

Lemma exists_lf_under : forall a s u v,
lf_beyond a s u v -> exists w,
(lf_beyond a s v w & lf_under a s u w).

Definition lf_equiv a s u v :=
(exists w, (lf_beyond a s u w & lf_beyond a s v w)).

Lemma lf_equiv_symm : forall a s u v,
lf_equiv a s u v -> lf_equiv a s v u.

Lemma lf_extend_id : forall a s u,
is_lf_symbol a s u ->
lf_extend a s (id a (lf_vertex u)) u = u.

Lemma lf_equiv_refl : forall a s u,
has_left_fractions a s ->
is_lf_symbol a s u -> lf_equiv a s u u.

Lemma lf_beyond_refl : forall a s u,
has_left_fractions a s ->
is_lf_symbol a s u -> lf_beyond a s u u.

Lemma lf_equiv_beyond_under : forall a s u v,
lf_equiv a s u v -> exists w,
(lf_beyond a s u w & lf_under a s v w).

Lemma exists_lf_further : forall a s u v w,
lf_beyond a s u v -> lf_under a s u w ->
exists z, (lf_beyond a s v z & lf_beyond a s w z).

Lemma lf_equiv_trans : forall a s u w,
has_left_fractions a s ->
(exists v, (lf_equiv a s u v & lf_equiv a s v w)) ->
lf_equiv a s u w.


Definition fills_in a s u v w :=
has_left_fractions a s &
is_lf_symbol a s u & is_lf_symbol a s v & is_lf_symbol a s w
& source u = target v & source w = lf_vertex v &
target w = lf_vertex u &
comp a (lf_forward w) (lf_backward v) =
comp a (lf_backward w) (lf_forward u).

Definition lf_filler a s u v :=
lf_choice a s (lf_backward v) (lf_forward u).

Definition lf_make_comp a (s:E) u v w :=
lf_symbol (comp a (lf_forward w) (lf_forward v))
(comp a (lf_backward w) (lf_backward u)).

Definition lf_comp_rep a s u v :=
lf_make_comp a s u v (lf_filler a s u v).

Lemma fills_in_lf_filler : forall a s u v,
has_left_fractions a s ->
is_lf_symbol a s u -> is_lf_symbol a s v ->
source u = target v ->
fills_in a s u v (lf_filler a s u v).

Lemma source_lf_make_comp : forall a s u v w,
fills_in a s u v w ->
source (lf_make_comp a s u v w) = source v.

Lemma target_lf_make_comp : forall a s u v w,
fills_in a s u v w ->
target (lf_make_comp a s u v w) = target u.

Lemma lf_vertex_lf_make_comp : forall a s u v w,
fills_in a s u v w ->
lf_vertex (lf_make_comp a s u v w) =
lf_vertex w.

Lemma is_lf_symbol_lf_make_comp : forall a s u v w,
fills_in a s u v w ->
is_lf_symbol a s (lf_make_comp a s u v w).

Lemma lf_forward_lf_make_comp : forall a s u v w,
lf_forward (lf_make_comp a s u v w) =
comp a (lf_forward w) (lf_forward v).

Lemma lf_backward_lf_make_comp : forall a s u v w,
lf_backward (lf_make_comp a s u v w) =
comp a (lf_backward w) (lf_backward u).

Lemma lf_beyond_fills_in : forall a s u v w,
(exists y, (fills_in a s u v y &
lf_beyond a s y w)) ->
fills_in a s u v w.

Lemma lf_beyond_lf_make_comp : forall a s u v w y,
fills_in a s u v w ->
lf_beyond a s w y ->
lf_beyond a s (lf_make_comp a s u v w) (lf_make_comp a s u v y).

Definition lf_lean_to a s e f g h i j :=
has_left_fractions a s &
mor a e & mor a f & mor a g & mor a h & mor a i & mor a j &
inc e s & inc h s &
source e = source f & source g = target e & source h = target f
& target g = target h & source i = target e & source j = target f
& target i = target j &
comp a g e = comp a h f & comp a i e = comp a j f.

Lemma show_lf_lean_to : forall a s e f g h i j,
has_left_fractions a s ->
mor a f -> mor a g -> mor a i -> mor a j ->
inc e s -> inc h s ->
source e = source f -> source g = target e -> source h = target f
-> source i = target e -> source j = target f ->
comp a g e = comp a h f -> comp a i e = comp a j f->
lf_lean_to a s e f g h i j.

Definition closes_lf_lean_to a s e f g h i j k l :=
lf_lean_to a s e f g h i j &
mor a k & mor a l & inc l s &
source k = target g & source l = target i &
target k = target l &
comp a k g = comp a l i & comp a k h = comp a l j.

Lemma lf_lean_to_closure : forall a s e f g h i j,
lf_lean_to a s e f g h i j ->
(exists k, exists l,(closes_lf_lean_to a s e f g h i j k l)).


Definition weakly_reps_lf_comp a s u v w :=
has_left_fractions a s &
is_lf_symbol a s u & is_lf_symbol a s v & is_lf_symbol a s w
& source u = target v & source w = source v &
target w = target u &
(exists p, exists q,( mor a p & mor a q &
source p = lf_vertex v & source q = lf_vertex u &
target p = lf_vertex w & target q = lf_vertex w &
comp a p (lf_backward v) = comp a q (lf_forward u) &
comp a p (lf_forward v) = lf_forward w &
comp a q (lf_backward u) = lf_backward w)).

Lemma weak_rep_lf_equiv : forall a s u v w y,
fills_in a s u v y ->
weakly_reps_lf_comp a s u v w ->
lf_equiv a s w (lf_make_comp a s u v y).

Lemma lf_beyond_weakly_reps_make_comp : forall a s u v w y,
lf_beyond a s u w -> fills_in a s w v y ->
weakly_reps_lf_comp a s u v (lf_make_comp a s w v y).

Lemma weakly_reps_make_comp : forall a s u v y,
fills_in a s u v y ->
weakly_reps_lf_comp a s u v (lf_make_comp a s u v y).

Lemma source_lf_comp_rep : forall a s u v,
has_left_fractions a s ->
is_lf_symbol a s u -> is_lf_symbol a s v ->
source u = target v ->
source (lf_comp_rep a s u v) = source v.

Lemma target_lf_comp_rep : forall a s u v,
has_left_fractions a s ->
is_lf_symbol a s u -> is_lf_symbol a s v ->
source u = target v ->
target (lf_comp_rep a s u v) = target u.

Lemma is_lf_symbol_lf_comp_rep : forall a s u v,
has_left_fractions a s ->
is_lf_symbol a s u -> is_lf_symbol a s v ->
source u = target v ->
is_lf_symbol a s (lf_comp_rep a s u v).

Lemma weakly_reps_lf_comp_lf_comp_rep : forall a s u v,
has_left_fractions a s -> is_lf_symbol a s u ->
is_lf_symbol a s v -> source u = target v ->
weakly_reps_lf_comp a s u v (lf_comp_rep a s u v).
 

Lemma lf_equiv_lf_comp_rep : forall a s u v w,
weakly_reps_lf_comp a s u v w ->
lf_equiv a s w (lf_comp_rep a s u v).

Lemma weakly_reps_lf_comp_equiv : forall a s w1 w2,
(exists u, exists v, (weakly_reps_lf_comp a s u v w1 &
weakly_reps_lf_comp a s u v w2)) ->
lf_equiv a s w1 w2.

Lemma weakly_reps_make_comp_beyond : forall a s u v u1 v1 y,
has_left_fractions a s -> lf_beyond a s u u1 ->
lf_beyond a s v v1 -> source u = target v ->
fills_in a s u1 v1 y ->
weakly_reps_lf_comp a s u v (lf_make_comp a s u1 v1 y).

Lemma weakly_reps_comp_rep_beyond : forall a s u v u1 v1,
has_left_fractions a s -> lf_beyond a s u u1 ->
lf_beyond a s v v1 -> source u = target v ->
weakly_reps_lf_comp a s u v (lf_comp_rep a s u1 v1).


Lemma lf_comp_indep : forall a s u v u1 v1,
has_left_fractions a s -> lf_equiv a s u u1 ->
lf_equiv a s v v1 -> source u = target v ->
lf_equiv a s (lf_comp_rep a s u v) (lf_comp_rep a s u1 v1).


Lemma lf_equiv_make_comp_comp_rep : forall a s u v y,
fills_in a s u v y ->
lf_equiv a s (lf_make_comp a s u v y) (lf_comp_rep a s u v).


Definition assoc_board a s u v w x y z :=
fills_in a s u v x &
fills_in a s v w y &
fills_in a s x y z.


Definition ffb_symbol a (s:E) y z :=
(lf_symbol (comp a (lf_forward z) (lf_forward y)) (lf_backward z)).


Definition fbb_symbol a (s:E) x z :=
(lf_symbol (lf_forward z) (comp a (lf_backward z) (lf_backward x))).
 
Definition ffb_symbol_facts (a s u v w x y z k:E) :=
lf_forward k = comp a (lf_forward z) (lf_forward y) &
lf_backward k = lf_backward z &
source k = lf_vertex w &
target k = lf_vertex x &
target (lf_forward k) = lf_vertex z &
target (lf_backward k) = lf_vertex z &
lf_symbol_like k &
mor a (lf_forward k) &
mor a (lf_backward k) &
inc (lf_backward k) s &
is_lf_symbol a s k &
fills_in a s (lf_make_comp a s u v x) w k.

Lemma get_ffb_symbol_facts : forall a s u v w x y z,
assoc_board a s u v w x y z ->
ffb_symbol_facts a s u v w x y z (ffb_symbol a s y z).

Definition fbb_symbol_facts (a s u v w x y z k:E) :=
lf_forward k = lf_forward z &
lf_backward k = comp a (lf_backward z) (lf_backward x) &
source k = lf_vertex y &
target k = lf_vertex u &
target (lf_forward k) = lf_vertex z &
target (lf_backward k) = lf_vertex z &
lf_symbol_like k &
mor a (lf_forward k) &
mor a (lf_backward k) &
inc (lf_backward k) s &
is_lf_symbol a s k &
fills_in a s u (lf_make_comp a s v w y) k.

Lemma get_fbb_symbol_facts : forall a s u v w x y z,
assoc_board a s u v w x y z ->
fbb_symbol_facts a s u v w x y z (fbb_symbol a s x z).

Lemma make_comp_assoc_board : forall a s u v w x y z,
assoc_board a s u v w x y z ->
lf_make_comp a s (lf_make_comp a s u v x) w (ffb_symbol a s y z) =
lf_make_comp a s u (lf_make_comp a s v w y) (fbb_symbol a s x z).

Lemma lf_comp_rep_assoc : forall a s u v w,
has_left_fractions a s ->
is_lf_symbol a s u -> is_lf_symbol a s v ->
is_lf_symbol a s w ->
source u = target v -> source v = target w ->
lf_equiv a s (lf_comp_rep a s (lf_comp_rep a s u v) w)
(lf_comp_rep a s u (lf_comp_rep a s v w)).

Definition left_id_filler a (s:E) u :=
lf_symbol (id a (lf_vertex u)) (lf_backward u).

Definition right_id_filler a (s:E) u :=
lf_symbol (lf_forward u) (id a (lf_vertex u)).

Lemma ob_lf_vertex : forall a u,
(exists s,(is_lf_symbol a s u)) ->
ob a (lf_vertex u).

Lemma source_left_id_filler : forall a s u,
is_lf_symbol a s u ->
source (left_id_filler a s u) = lf_vertex u.

Lemma target_left_id_filler : forall a s u,
is_lf_symbol a s u ->
target (left_id_filler a s u) = target u.

Lemma source_right_id_filler : forall a s u,
is_lf_symbol a s u ->
source (right_id_filler a s u) = source u.

Lemma target_right_id_filler : forall a s u,
is_lf_symbol a s u ->
target (right_id_filler a s u) = lf_vertex u.

Lemma lf_forward_left_id_filler : forall a s u,
is_lf_symbol a s u ->
lf_forward (left_id_filler a s u) = id a (lf_vertex u).

Lemma lf_backward_left_id_filler : forall a s u,
is_lf_symbol a s u ->
lf_backward (left_id_filler a s u) = lf_backward u.

Lemma lf_forward_right_id_filler : forall a s u,
is_lf_symbol a s u ->
lf_forward (right_id_filler a s u) = lf_forward u.

Lemma lf_backward_right_id_filler : forall a s u,
is_lf_symbol a s u ->
lf_backward (right_id_filler a s u) = id a (lf_vertex u).

Lemma lf_vertex_left_id_filler : forall a s u,
is_lf_symbol a s u ->
lf_vertex (left_id_filler a s u) = lf_vertex u.

Lemma lf_vertex_right_id_filler : forall a s u,
is_lf_symbol a s u ->
lf_vertex (right_id_filler a s u) = lf_vertex u.

Lemma is_lf_symbol_left_id_filler : forall a s u,
is_lf_symbol a s u ->
is_lf_symbol a s (left_id_filler a s u).

Lemma is_lf_symbol_right_id_filler : forall a s u,
has_left_fractions a s -> is_lf_symbol a s u ->
is_lf_symbol a s (right_id_filler a s u).

Lemma lf_forward_lf_id_rep : forall a x,
lf_forward (lf_id_rep a x) = id a x.

Lemma lf_backward_lf_id_rep : forall a x,
lf_backward (lf_id_rep a x) = id a x.

Lemma fills_in_left_id_filler : forall a s u,
has_left_fractions a s -> is_lf_symbol a s u ->
fills_in a s (lf_id_rep a (target u)) u (left_id_filler a s u).

Lemma fills_in_right_id_filler : forall a s u,
has_left_fractions a s -> is_lf_symbol a s u ->
fills_in a s u (lf_id_rep a (source u)) (right_id_filler a s u).

Lemma lf_make_comp_left_id_filler : forall a s u,
is_lf_symbol a s u ->
lf_make_comp a s (lf_id_rep a (target u)) u
(left_id_filler a s u) = u.

Lemma lf_make_comp_right_id_filler : forall a s u,
is_lf_symbol a s u ->
lf_make_comp a s u (lf_id_rep a (source u))
(right_id_filler a s u) = u.

Lemma lf_left_id : forall a s u,
has_left_fractions a s ->
is_lf_symbol a s u ->
lf_equiv a s (lf_comp_rep a s (lf_id_rep a (target u)) u) u.

Lemma lf_right_id : forall a s u,
has_left_fractions a s ->
is_lf_symbol a s u ->
lf_equiv a s (lf_comp_rep a s u (lf_id_rep a (source u))) u.


Definition taut_lf_symbol a y :=
lf_symbol y (id a (target y)).

Lemma source_taut_lf_symbol : forall a y,
source (taut_lf_symbol a y) = source y.

Lemma target_taut_lf_symbol : forall a y,
mor a y ->
target (taut_lf_symbol a y) = target y.

Lemma lf_forward_taut_lf_symbol : forall a y,
lf_forward (taut_lf_symbol a y) = y.

Lemma lf_backward_taut_lf_symbol : forall a y,
lf_backward (taut_lf_symbol a y) = id a (target y).

Lemma lf_vertex_taut_lf_symbol : forall a y,
mor a y ->
lf_vertex (taut_lf_symbol a y) = target y.

Lemma is_lf_symbol_taut_lf_symbol : forall a s y,
has_left_fractions a s -> mor a y ->
is_lf_symbol a s (taut_lf_symbol a y).

Lemma taut_lf_symbol_id : forall a x,
ob a x ->
taut_lf_symbol a (id a x) = lf_id_rep a x.

Lemma fills_in_taut_lf_symbol : forall a s y z,
has_left_fractions a s ->
mor a y -> mor a z -> source y = target z ->
fills_in a s (taut_lf_symbol a y)
(taut_lf_symbol a z) (taut_lf_symbol a y).

Lemma comp_taut_lf_symbol : forall a s y z,
has_left_fractions a s ->
mor a y -> mor a z -> source y = target z ->
lf_equiv a s
(lf_comp_rep a s (taut_lf_symbol a y) (taut_lf_symbol a z))
(taut_lf_symbol a (comp a y z)).



Definition inverse_lf_symbol a y :=
lf_symbol (id a (target y)) y.

Lemma source_inverse_lf_symbol : forall a y,
mor a y ->
source (inverse_lf_symbol a y) = target y.

Lemma target_inverse_lf_symbol : forall a y,
mor a y ->
target (inverse_lf_symbol a y) = source y.

Lemma lf_forward_inverse_lf_symbol : forall a y,
lf_forward (inverse_lf_symbol a y) = id a (target y).

Lemma lf_backward_inverse_lf_symbol : forall a y,
lf_backward (inverse_lf_symbol a y) = y.

Lemma lf_vertex_inverse_lf_symbol : forall a y,
lf_vertex (inverse_lf_symbol a y) = target y.

Lemma is_lf_symbol_inverse_lf_symbol : forall a s y,
has_left_fractions a s -> inc y s ->
is_lf_symbol a s (inverse_lf_symbol a y).

Lemma fills_in_left_lf_inverse : forall a s y,
has_left_fractions a s -> inc y s ->
fills_in a s (inverse_lf_symbol a y) (taut_lf_symbol a y)
(lf_id_rep a (target y)).

Lemma fills_in_right_lf_inverse : forall a s y,
has_left_fractions a s -> inc y s ->
fills_in a s (taut_lf_symbol a y) (inverse_lf_symbol a y)
(lf_id_rep a (target y)).

Lemma lf_left_inverse : forall a s y,
has_left_fractions a s -> inc y s ->
lf_equiv a s
(lf_comp_rep a s (inverse_lf_symbol a y) (taut_lf_symbol a y))
(lf_id_rep a (source y)).

Lemma lf_right_inverse : forall a s y,
has_left_fractions a s -> inc y s ->
lf_equiv a s
(lf_comp_rep a s (taut_lf_symbol a y) (inverse_lf_symbol a y))
(lf_id_rep a (target y)).

Lemma lf_symbol_equiv : forall a s u,
has_left_fractions a s -> is_lf_symbol a s u ->
lf_equiv a s u
(lf_comp_rep a s
(inverse_lf_symbol a (lf_backward u))
(taut_lf_symbol a (lf_forward u))).

End Left_Fractions.

\end{verbatim}

\subsection{The module {\tt Left\_Fraction\_Category}}

\begin{verbatim}
Module Left_Fraction_Category.
Export Left_Fractions.
Export Associating_Quotient.

Definition lf_symbol_container a s :=
Image.create (Cartesian.product (morphisms a) s)
(fun z => lf_symbol (pr1 z) (pr2 z)).

Lemma inc_lf_symbol_container : forall a s u,
is_lf_symbol a s u ->
inc u (lf_symbol_container a s).

Definition lf_symbol_set a s :=
Z (lf_symbol_container a s) (is_lf_symbol a s).

Lemma inc_lf_symbol_set : forall a s u,
inc u (lf_symbol_set a s) = is_lf_symbol a s u.

Definition lfc_rqcat a s :=
Category.Notations.create (objects a)
(lf_symbol_set a s) (lf_comp_rep a s) (lf_id_rep a) (structure a).

Lemma is_ob_lfc_rqcat : forall a s x,
has_left_fractions a s ->
is_ob (lfc_rqcat a s) x = ob a x.

Lemma is_mor_lfc_rqcat : forall a s u,
has_left_fractions a s ->
is_mor (lfc_rqcat a s) u = is_lf_symbol a s u.

Lemma comp_lfc_rqcat : forall a s u v,
has_left_fractions a s ->
is_lf_symbol a s u ->
is_lf_symbol a s v ->
source u = target v ->
comp (lfc_rqcat a s) u v = lf_comp_rep a s u v.

Lemma id_lfc_rqcat : forall a s x,
has_left_fractions a s -> ob a x ->
id (lfc_rqcat a s) x = lf_id_rep a x.

Lemma rqcat_lfc_rqcat : forall a s,
has_left_fractions a s ->
rqcat (lfc_rqcat a s).

Definition lfer a s :=
Z (Cartesian.product (lf_symbol_set a s) (lf_symbol_set a s))
(fun z => lf_equiv a s (pr1 z) (pr2 z)).

Lemma related_lfer : forall a s u v,
has_left_fractions a s ->
related (lfer a s) u v = lf_equiv a s u v.

Lemma related_lfer_rw : forall a s u v,
has_left_fractions a s ->
related (lfer a s) u v =
(is_lf_symbol a s u & is_lf_symbol a s v &
source u = source v & target u = target v &
lf_equiv a s u v).

Lemma inc_lfer : forall a s x,
has_left_fractions a s ->
inc x (lfer a s) =
(is_pair x & (related (lfer a s) (pr1 x) (pr2 x))).

Lemma lf_equiv_properties : forall a s u v,
lf_equiv a s u v ->
(lf_equiv a s u v &
is_lf_symbol a s u & is_lf_symbol a s v &
ob a (source u) & ob a (target u) &
ob a (source v) & ob a (target v) &
source u = source v & target u = target v).

Lemma lf_ob_target : forall a u,
(exists s, (is_lf_symbol a s u)) ->
ob a (target u).

Lemma lf_ob_source : forall a u,
(exists s, (is_lf_symbol a s u)) ->
ob a (source u).

Lemma rqcat_equiv_rel_lfer : forall a s,
has_left_fractions a s ->
rqcat_equiv_rel (lfc_rqcat a s) (lfer a s).

Definition left_frac_cat a s :=
quotient_cat (lfc_rqcat a s) (lfer a s).

Lemma left_frac_cat_axioms : forall a s,
has_left_fractions a s ->
Category.axioms (left_frac_cat a s).

Definition lf_class a s u := arrow_class (lfer a s) u.

Lemma source_lf_class : forall a s u,
source (lf_class a s u) = source u.

Lemma target_lf_class : forall a s u,
target (lf_class a s u) = target u.

Lemma ob_left_frac_cat : forall a s x,
has_left_fractions a s ->
ob (left_frac_cat a s) x = ob a x.

Lemma mor_left_frac_cat : forall a s v,
has_left_fractions a s ->
mor (left_frac_cat a s) v =
(exists u, (is_lf_symbol a s u & v = lf_class a s u)).

Lemma comp_left_frac_cat_lf_class : forall a s u v,
has_left_fractions a s ->
is_lf_symbol a s u -> is_lf_symbol a s v ->
source u = target v ->
comp (left_frac_cat a s) (lf_class a s u) (lf_class a s v) =
lf_class a s (lf_comp_rep a s u v).

Lemma id_left_frac_cat : forall a s x,
has_left_fractions a s -> ob a x ->
id (left_frac_cat a s) x = lf_class a s (lf_id_rep a x).

Lemma eq_lf_class : forall a s u v,
has_left_fractions a s ->
is_lf_symbol a s u -> is_lf_symbol a s v ->
(lf_class a s u = lf_class a s v) = (lf_equiv a s u v).

Definition lf_taut a s u :=
lf_class a s (taut_lf_symbol a u).

Definition lf_inverse a s u :=
lf_class a s (inverse_lf_symbol a u).

Lemma source_lf_taut : forall a s u,
source (lf_taut a s u) = source u.

Lemma target_lf_taut : forall a s u,
mor a u ->
target (lf_taut a s u) = target u.

Lemma source_lf_inverse : forall a s u,
mor a u ->
source (lf_inverse a s u) = target u.

Lemma target_lf_inverse : forall a s u,
mor a u ->
target (lf_inverse a s u) = source u.

Lemma mor_lf_taut : forall a s u,
has_left_fractions a s -> mor a u ->
mor (left_frac_cat a s) (lf_taut a s u).

Lemma mor_lf_inverse : forall a s u,
has_left_fractions a s -> inc u s ->
mor (left_frac_cat a s) (lf_inverse a s u).

Lemma lf_taut_id : forall a s x,
has_left_fractions a s -> ob a x ->
lf_taut a s (id a x) = id (left_frac_cat a s) x.

Lemma comp_lf_class : forall a s u v,
has_left_fractions a s ->
is_lf_symbol a s u -> is_lf_symbol a s v ->
source u = target v ->
comp (left_frac_cat a s) (lf_class a s u) (lf_class a s v) =
lf_class a s (lf_comp_rep a s u v).

Lemma comp_lf_taut : forall a s u v,
has_left_fractions a s ->
mor a u -> mor a v -> source u = target v ->
comp (left_frac_cat a s) (lf_taut a s u) (lf_taut a s v) =
lf_taut a s (comp a u v).

Lemma are_inverse_lf_taut_lf_inverse : forall a s u,
has_left_fractions a s -> inc u s ->
are_inverse (left_frac_cat a s) (lf_taut a s u) (lf_inverse a s u).

Lemma is_lf_symbol_arrow_rep : forall a s y,
has_left_fractions a s ->
mor (left_frac_cat a s) y ->
is_lf_symbol a s (arrow_rep y).

Lemma lf_class_arrow_rep : forall a s y,
has_left_fractions a s ->
mor (left_frac_cat a s) y ->
lf_class a s (arrow_rep y) = y.

Lemma lf_source_arrow_rep : forall y,
(exists a, exists s, (has_left_fractions a s &
mor (left_frac_cat a s) y)) ->
source (arrow_rep y) = source y.

Lemma lf_target_arrow_rep : forall y,
(exists a, exists s, (has_left_fractions a s &
mor (left_frac_cat a s) y)) ->
target (arrow_rep y) = target y.

Lemma comp_etc_arrow_rep : forall a s y,
has_left_fractions a s ->
mor (left_frac_cat a s) y ->
comp (left_frac_cat a s)
(lf_inverse a s (lf_backward (arrow_rep y)))
(lf_taut a s (lf_forward (arrow_rep y))) = y.

Lemma comp_etc_lf_symbol : forall a s k,
has_left_fractions a s ->
is_lf_symbol a s k ->
comp (left_frac_cat a s)
(lf_inverse a s (lf_backward k))
(lf_taut a s (lf_forward k)) = lf_class a s k.

Lemma lfc_mor_expression : forall a s y,
has_left_fractions a s ->
mor (left_frac_cat a s) y ->
(exists u, exists v, (mor a u & inc v s & target u = target v &
y = comp (left_frac_cat a s) (lf_inverse a s v) (lf_taut a s u))).


Definition lf_proj a s :=
Functor.create a (left_frac_cat a s) (lf_taut a s).

Lemma source_lf_proj : forall a s,
source (lf_proj a s) = a.

Lemma target_lf_proj : forall a s,
target (lf_proj a s) = (left_frac_cat a s).

Lemma lf_proj_property : forall a s,
has_left_fractions a s ->
Functor.property a (left_frac_cat a s) (fun (x:E)=> x)
(lf_taut a s).

Lemma lf_proj_axioms : forall a s,
has_left_fractions a s ->
Functor.axioms (lf_proj a s).

Lemma fob_lf_proj : forall a s x,
has_left_fractions a s -> ob a x ->
fob (lf_proj a s) x = x.

Lemma fmor_lf_proj : forall a s y,
has_left_fractions a s -> mor a y ->
fmor (lf_proj a s) y = lf_taut a s y.

Lemma invertible_fmor_lf_proj : forall a s y,
has_left_fractions a s -> inc y s ->
invertible (left_frac_cat a s) (fmor (lf_proj a s) y).

Lemma inverse_fmor_lf_proj : forall a s y,
has_left_fractions a s -> inc y s ->
inverse (left_frac_cat a s) (fmor (lf_proj a s) y) =
lf_inverse a s y.


Definition lf_dotted_situation a s f :=
has_left_fractions a s &
Functor.axioms f &
source f = a &
(forall y, inc y s -> invertible (target f) (fmor f y)).

Definition lf_symb_dot_situation a s f u :=
lf_dotted_situation a s f &
is_lf_symbol a s u.

Definition lf_symb_dot_facts a s f u:=
lf_symb_dot_situation a s f u &
is_lf_symbol a s u &
lf_dotted_situation a s f &
has_left_fractions a s &
Functor.axioms f &
source f = a &
(forall y, inc y s -> invertible (target f) (fmor f y)) &
mor (source f) (lf_forward u) &
mor (source f) (lf_backward u) &
inc (lf_backward u) s &
mor (target f) (fmor f (lf_forward u)) &
mor (target f) (fmor f (lf_backward u)) &
lf_symbol_like u &
source (fmor f (lf_forward u)) = fob f (source u) &
source (fmor f (lf_backward u)) = fob f (target u) &
target (fmor f (lf_forward u)) = fob f (lf_vertex u) &
target (fmor f (lf_backward u)) = fob f (lf_vertex u) &
invertible (target f) (fmor f (lf_backward u)) &
mor (target f) (inverse (target f) (fmor f (lf_backward u))) &
source (inverse (target f) (fmor f (lf_backward u))) =
fob f (lf_vertex u).

Lemma lf_symb_dot_situation_rw : forall a s f u,
lf_symb_dot_situation a s f u = lf_symb_dot_facts a s f u.

Definition lf_symb_dot (a s:E) f u :=
comp (target f)
(inverse (target f) (fmor f (lf_backward u)))
(fmor f (lf_forward u)).

Lemma source_lf_symb_dot : forall a s f u,
lf_symb_dot_situation a s f u ->
source (lf_symb_dot a s f u) = fob f (source u).

Lemma target_lf_symb_dot : forall a s f u,
lf_symb_dot_situation a s f u ->
target (lf_symb_dot a s f u) = fob f (target u).

Lemma mor_lf_symb_dot : forall a s f u,
lf_symb_dot_situation a s f u ->
mor (target f) (lf_symb_dot a s f u).

Lemma lf_symb_dot_beyond_invariant : forall a s f u v,
lf_symb_dot_situation a s f u ->
lf_beyond a s u v ->
lf_symb_dot a s f u = lf_symb_dot a s f v.

Lemma lf_symb_dot_equiv_invariant : forall a s f u v,
lf_symb_dot_situation a s f u ->
lf_equiv a s u v ->
lf_symb_dot a s f u = lf_symb_dot a s f v.

Lemma lf_symb_dot_make_comp : forall a s f u v w,
lf_symb_dot_situation a s f u ->
lf_symb_dot_situation a s f v ->
fills_in a s u v w ->
lf_symb_dot a s f (lf_make_comp a s u v w) =
comp (target f) (lf_symb_dot a s f u) (lf_symb_dot a s f v).

Lemma lf_symb_dot_comp_rep : forall a s f u v,
lf_symb_dot_situation a s f u ->
lf_symb_dot_situation a s f v ->
source u = target v ->
lf_symb_dot a s f (lf_comp_rep a s u v) =
comp (target f) (lf_symb_dot a s f u) (lf_symb_dot a s f v).

Lemma lf_symb_dot_id_rep : forall a s f x,
lf_dotted_situation a s f ->
ob a x ->
lf_symb_dot a s f (lf_id_rep a x) = id (target f) (fob f x).

Definition lf_dotted a s f :=
Functor.create (left_frac_cat a s) (target f)
(fun u => lf_symb_dot a s f (arrow_rep u)).

Lemma source_lf_dotted : forall a s f,
source (lf_dotted a s f) = left_frac_cat a s.

Lemma target_lf_dotted : forall a s f,
target (lf_dotted a s f) = target f.

Lemma lf_equiv_arrow_rep_lf_class : forall a s u,
has_left_fractions a s ->
is_lf_symbol a s u ->
lf_equiv a s u (arrow_rep (lf_class a s u)).

Lemma lf_dotted_property : forall a s f,
lf_dotted_situation a s f ->
Functor.property (left_frac_cat a s) (target f)
(fob f) (fun u => lf_symb_dot a s f (arrow_rep u)).

Lemma lf_dotted_axioms : forall a s f,
lf_dotted_situation a s f ->
Functor.axioms (lf_dotted a s f).

Lemma fob_lf_dotted : forall a s f x,
lf_dotted_situation a s f ->
ob a x ->
fob (lf_dotted a s f) x = fob f x.

Lemma fmor_lf_dotted : forall a s f u,
mor (left_frac_cat a s) u ->
fmor (lf_dotted a s f) u = lf_symb_dot a s f (arrow_rep u).

Lemma fmor_lf_dotted_lf_class : forall a s f u,
lf_dotted_situation a s f ->
is_lf_symbol a s u ->
fmor (lf_dotted a s f) (lf_class a s u) = lf_symb_dot a s f u.

Lemma fmor_lf_dotted_lf_taut : forall a s f y,
lf_dotted_situation a s f ->
mor a y ->
fmor (lf_dotted a s f) (lf_taut a s y) = fmor f y.

Lemma fmor_lf_dotted_lf_inverse : forall a s f y,
lf_dotted_situation a s f ->
inc y s ->
fmor (lf_dotted a s f) (lf_inverse a s y) =
inverse (target f) (fmor f y).


Lemma fcompose_lf_dotted_lf_proj : forall a s f,
lf_dotted_situation a s f ->
fcompose (lf_dotted a s f) (lf_proj a s) = f.


Lemma lf_dotted_like_unique : forall a s g h,
has_left_fractions a s ->
Functor.axioms g -> Functor.axioms h ->
source g = left_frac_cat a s ->
source h = left_frac_cat a s ->
(fcompose g (lf_proj a s) = fcompose h (lf_proj a s)) ->
g = h.

End Left_Fraction_Category.

\end{verbatim}

\newpage

\section{The file {\tt gzloc.v}}

\begin{verbatim}

Require Export left_fractions.

\end{verbatim}

\subsection{The module {\tt GZ\_Localization}}

\begin{verbatim}
Module GZ_Localization.
Export GZ_Thm.
Export Left_Fraction_Category.

Lemma localizing_system_recall : forall a s,
localizing_system a s =
(Category.axioms a &
(forall u, inc u s -> mor a u)).

Lemma multiplicative_system_recall : forall a s,
multiplicative_system a s =
(localizing_system a s &
(forall y z, inc y s -> inc z s ->
source y = target z -> inc (comp a y z) s)).

Definition localizes a s f :=
localizing_system a s &
Functor.axioms f &
source f = a &
(forall y, inc y s -> invertible (target f) (fmor f y)).

Definition completes_triangle f g h :=
Functor.axioms f & Functor.axioms g & Functor.axioms h &
source f = source g & source h = target f &
target h = target g & fcompose h f = g.

Definition dotted_choice f g :=
choose (completes_triangle f g).

Lemma completes_triangle_dotted_choice : forall f g,
completes_triangle f g (dotted_choice f g) =
(exists h, (completes_triangle f g h)).

Definition is_localization a s f :=
localizes a s f &
(forall g, (localizes a s g ->
completes_triangle f g (dotted_choice f g))) &
(forall g h, completes_triangle f g h ->
h = dotted_choice f g).

Lemma dotted_choice_axioms : forall f g,
(exists a, exists s, (is_localization a s f & localizes a s g)) ->
Functor.axioms (dotted_choice f g).

Lemma source_dotted_choice : forall f g,
(exists a, exists s, (is_localization a s f & localizes a s g)) ->
source (dotted_choice f g) = target f.

Lemma target_dotted_choice : forall f g,
(exists a, exists s, (is_localization a s f & localizes a s g)) ->
target (dotted_choice f g) = target g.

Lemma fcompose_dotted_choice : forall f g,
(exists a, exists s, (is_localization a s f & localizes a s g)) ->
fcompose (dotted_choice f g) f = g.

Lemma localizes_fcompose : forall a s f g,
localizes a s f -> Functor.axioms g -> source g = target f ->
localizes a s (fcompose g f).

Lemma eq_dotted_choice : forall f g h,
completes_triangle f g h ->
(exists a, exists s, (is_localization a s f)) ->
h = dotted_choice f g.

Lemma fcompose_dotted_choice_dotted_choice : forall f g h,
(exists a, exists s, (is_localization a s f & is_localization a s g &
localizes a s h)) ->
fcompose (dotted_choice g h) (dotted_choice f g) =
dotted_choice f h.

Lemma dotted_choice_refl : forall f,
(exists a, exists s, (is_localization a s f)) ->
dotted_choice f f = fidentity (target f).


Lemma are_finverse_dotted_choice : forall f g,
(exists a, exists s, 
(is_localization a s f & is_localization a s g)) ->
are_finverse (dotted_choice f g) (dotted_choice g f).

Lemma has_finverse_dotted_choice : forall f g,
(exists a, exists s, 
(is_localization a s f & is_localization a s g)) ->
has_finverse (dotted_choice f g).

Lemma finverse_dotted_choice : forall f g,
(exists a, exists s, 
(is_localization a s f & is_localization a s g)) ->
finverse (dotted_choice f g) = dotted_choice g f.

Lemma is_localization_gz_proj : forall a s,
localizing_system a s ->
is_localization a s (gz_proj a s).

Lemma is_localization_lf_proj : forall a s,
has_left_fractions a s ->
is_localization a s (lf_proj a s).

Lemma are_finverse_dotted_choice_gz_proj_lf_proj : forall a s,
has_left_fractions a s ->
are_finverse (dotted_choice (gz_proj a s) (lf_proj a s))
(dotted_choice (lf_proj a s) (gz_proj a s)).


Definition oppms s := Image.create s flip.

Lemma inc_oppms : forall s y,
inc y (oppms s) = inc (flip y) s.

Lemma oppms_oppms : forall s, oppms (oppms s) = s.

Lemma localizing_system_oppms : forall a s,
localizing_system (opp a) (oppms s) = localizing_system a s.

Lemma multiplicative_system_oppms : forall a s,
multiplicative_system (opp a) (oppms s) =
multiplicative_system a s.

Lemma are_inverse_opp : forall a u v,
are_inverse (opp a) (flip u) (flip v) =
are_inverse a u v.

Lemma invertible_opp : forall a u,
invertible (opp a) (flip u) =
invertible a u.

Lemma localizes_oppms_oppf : forall a s f,
localizes (opp a) (oppms s) (oppf f) =
localizes a s f.

Lemma completes_triangle_oppf : forall f g h,
completes_triangle (oppf f) (oppf g) (oppf h) =
completes_triangle f g h.

Lemma is_localization_oppms_oppf : forall a s f,
is_localization (opp a) (oppms s) (oppf f) =
is_localization a s f.


Lemma has_left_fractions_rw : forall a s,
has_left_fractions a s = (
multiplicative_system a s
&
(forall x, ob a x -> inc (id a x) s)
&
(forall r g, inc r s -> mor a g -> source r = source g ->
exists p, exists q, (mor a p & inc q s
& target p = target q & source p = target r
& source q = target g & comp a q g =
comp a p r))
&
(forall v r t, inc v s -> mor a r -> mor a t ->
source r = target v -> source t = target v ->
comp a r v = comp a t v ->
exists w, (inc w s & source w = target r & source w = target t
& comp a w r = comp a w t))).

Definition has_right_fractions a s :=
multiplicative_system a s
&
(forall x, ob a x -> inc (id a x) s)
&
(forall r g, inc r s -> mor a g -> target r = target g ->
exists p, exists q, (mor a q & inc p s &
source p = source q & target p = source g &
target q = source r & comp a g p =
comp a r q))
&
(forall v r t, inc v s -> mor a r -> mor a t ->
source v = target r -> source v = target t ->
comp a v r = comp a v t ->
exists w, (inc w s & source r = target w & source t = target w
& comp a r w = comp a t w)).

Lemma has_right_fractions_rw : forall a s,
has_right_fractions a s =
has_left_fractions (opp a) (oppms s).

Definition right_frac_cat a s :=
opp (left_frac_cat (opp a) (oppms s)).

Definition rf_proj a s :=
oppf (lf_proj (opp a) (oppms s)).

Lemma right_frac_cat_axioms : forall a s,
has_right_fractions a s ->
Category.axioms (right_frac_cat a s).

Lemma rf_proj_axioms : forall a s,
has_right_fractions a s ->
Functor.axioms (rf_proj a s).

Lemma source_rf_proj : forall a s,
has_right_fractions a s ->
source (rf_proj a s) = a.

Lemma target_rf_proj : forall a s,
has_right_fractions a s ->
target (rf_proj a s) = right_frac_cat a s.

Lemma is_localization_rf_proj : forall a s,
has_right_fractions a s ->
is_localization a s (rf_proj a s).


Definition lf_vee a s p q :=
localizing_system a s &
mor a p & inc q s & target p = target q.

Definition lf_vee_image f p q :=
comp (target f) (inverse (target f) (fmor f q))
(fmor f p).


Definition lf_vee_equivalent a s p q r t :=
lf_vee a s p q & lf_vee a s r t &
source p = source r & source q = source t &
(exists y, exists z, (mor a y & mor a z &
source y = target p & source z = target r &
target y = target z &
comp a y p = comp a z r &
comp a y q = comp a z t &
inc (comp a y q) s)).

Definition left_fraction_description a s f :=
localizes a s f &
ob_iso f
&
(forall y, mor (target f) y ->
exists p, exists q, (lf_vee a s p q &
y = lf_vee_image f p q))
&
(forall p q r t, lf_vee a s p q -> lf_vee a s r t ->
(lf_vee_image f p q = lf_vee_image f r t) ->
lf_vee_equivalent a s p q r t).


Definition rf_wedge a s p q :=
localizing_system a s &
mor a p & inc q s & source p = source q.

Definition rf_wedge_image f p q :=
comp (target f) (fmor f p) (inverse (target f) (fmor f q)).

Definition rf_wedge_equivalent a s p q r t :=
rf_wedge a s p q & rf_wedge a s r t &
target p = target r & target q = target t &
(exists y, exists z, (mor a y & mor a z &
target y = source p & target z = source r &
source y = source z &
comp a p y = comp a r z &
comp a q y = comp a t z &
inc (comp a q y) s)).

Definition right_fraction_description a s f :=
localizes a s f &
ob_iso f
&
(forall y, mor (target f) y ->
exists p, exists q, (rf_wedge a s p q &
y = rf_wedge_image f p q))
&
(forall p q r t, rf_wedge a s p q -> rf_wedge a s r t ->
(rf_wedge_image f p q = rf_wedge_image f r t) ->
rf_wedge_equivalent a s p q r t).


Lemma lf_vee_image_lf_proj : forall a s p q,
has_left_fractions a s ->
lf_vee a s p q ->
lf_vee_image (lf_proj a s) p q =
lf_class a s (lf_symbol p q).

Lemma left_fraction_description_lf_proj : forall a s,
has_left_fractions a s ->
left_fraction_description a s (lf_proj a s).

Lemma mor_lf_vee_image : forall f p q,
(exists a, exists s, (lf_vee a s p q &
localizes a s f)) ->
mor (target f) (lf_vee_image f p q).

Lemma lf_vee_image_fcompose : forall f g p q,
Functor.axioms f -> Functor.axioms g -> source f = target g ->
mor (source g) p -> mor (source g) q ->
invertible (source f) (fmor g q) ->
target p = target q ->
lf_vee_image (fcompose f g) p q =
fmor f (lf_vee_image g p q).

Lemma left_fraction_description_invariant : forall a s f g,
left_fraction_description a s f ->
has_finverse g ->
source g = target f ->
left_fraction_description a s (fcompose g f).

Lemma left_fraction_description_for_loc : forall a s f,
has_left_fractions a s ->
is_localization a s f ->
left_fraction_description a s f.

Lemma left_fraction_description_gz_proj :
forall a s,
has_left_fractions a s ->
left_fraction_description a s (gz_proj a s).

Lemma inverse_opp : forall a u,
invertible a (flip u) ->
inverse (opp a) u = flip (inverse a (flip u)).

Lemma right_fraction_description_rw : forall a s f,
right_fraction_description a s f =
left_fraction_description (opp a) (oppms s) (oppf f).

Lemma right_fraction_description_for_loc : forall a s f,
has_right_fractions a s ->
is_localization a s f ->
right_fraction_description a s f.

Lemma right_fraction_description_gz_proj : forall a s,
has_right_fractions a s ->
right_fraction_description a s (gz_proj a s).


End GZ_Localization. 

\end{verbatim}

\newpage

\section{The file {\tt lfcx.v}}

\begin{verbatim}

Require Export updateA.
Require Export gzloc.

\end{verbatim}

\subsection{The module {\tt Coarse\_Cat}}

\begin{verbatim}
Module Coarse_Cat.
Export GZ_Localization.

Definition null_arrow a b := Arrow.create a b emptyset.

Lemma source_null_arrow : forall a b, source (null_arrow a b) = a.

Lemma target_null_arrow : forall a b, target (null_arrow a b) = b.

Lemma arrow_null_arrow : forall a b,
arrow (null_arrow a b) = emptyset.

Definition na_comp u v := null_arrow (source v) (target u).

Lemma source_na_comp : forall u v, source (na_comp u v) = source v.

Lemma target_na_comp : forall u v, target (na_comp u v) = target u.

Definition is_null_arrow u := u = null_arrow (source u) (target u).

Lemma is_null_arrow_null_arrow : forall a b,
is_null_arrow (null_arrow a b).

Lemma is_null_arrow_na_comp : forall u v,
is_null_arrow (na_comp u v).

Lemma null_arrow_extensionality : forall a b,
is_null_arrow a -> is_null_arrow b ->
source a = source b -> target a = target b -> a = b.

Lemma na_comp_assoc : forall a b c,
na_comp (na_comp a b) c = na_comp a (na_comp b c).

Lemma na_comp_null_arrow : forall a b c d,
na_comp (null_arrow a b) (null_arrow c d) =
null_arrow c b.

Lemma na_left_id : forall a u,
is_null_arrow a ->
u = null_arrow (target a) (target a) ->
na_comp u a = a.

Lemma na_right_id : forall a u,
is_null_arrow a ->
u = null_arrow (source a) (source a) ->
na_comp a u = a.

 
Definition coarse_arrows z :=
Image.create (product z z)
(fun x => null_arrow (pr1 x) (pr2 x)).

Lemma inc_coarse_arrows : forall u z,
inc u (coarse_arrows z) =
(is_null_arrow u & inc (source u) z & inc (target u) z).

Definition coarse_cat z :=
Category.Notations.create z (coarse_arrows z)
na_comp (fun x => null_arrow x x) emptyset.

Lemma is_ob_coarse_cat : forall x z,
is_ob (coarse_cat z) x = inc x z.

Lemma is_mor_coarse_cat : forall u z,
is_mor (coarse_cat z) u =
(is_null_arrow u & inc (source u) z & inc (target u) z).

Lemma comp_coarse_cat1 : forall u v z,
is_mor (coarse_cat z) u -> is_mor (coarse_cat z) v ->
source u = target v ->
comp (coarse_cat z) u v = na_comp u v.

Lemma id_coarse_cat : forall x z,
inc x z -> id (coarse_cat z) x = null_arrow x x.

Lemma coarse_cat_axioms : forall z,
Category.axioms (coarse_cat z).

Lemma ob_coarse_cat : forall x z,
ob (coarse_cat z) x = inc x z.

Lemma mor_coarse_cat : forall u z,
mor (coarse_cat z) u =
(is_null_arrow u & inc (source u) z & inc (target u) z).

Lemma comp_coarse_cat : forall u v z,
mor (coarse_cat z) u -> mor (coarse_cat z) v ->
source u = target v ->
comp (coarse_cat z) u v = na_comp u v.

End Coarse_Cat.

\end{verbatim}

\subsection{The module {\tt Lf\_Counterexample}}

\begin{verbatim}
Module Lf_Counterexample.
Export Coarse_Cat.
Export Left_Fraction_Category.

Inductive cx_ob : E -> Prop :=
is0 : cx_ob (R 0) |
is1 : cx_ob (R 1) |
is2 : cx_ob (R 2).

Ltac is_cx_ob :=
match goal with
| |- cx_ob (R 0) => ap is0
| |- cx_ob (R 1) => ap is1
| |- cx_ob (R 2) => ap is2
| _ => fail end.

Lemma cx_ob_rw : forall x,
cx_ob x = inc x (R 3).

Notation na00 := (null_arrow (R 0) (R 0)).
Notation na11 := (null_arrow (R 1) (R 1)).
Notation na22 := (null_arrow (R 2) (R 2)).

Notation na01 := (null_arrow (R 0) (R 1)).
Notation na02 := (null_arrow (R 0) (R 2)).

Notation na12 := (null_arrow (R 1) (R 2)).
Notation na21 := (null_arrow (R 2) (R 1)).

Inductive cx_mor : E -> Prop :=
is00 : cx_mor na00 |
is11 : cx_mor na11 |
is22 : cx_mor na22 |
is01 : cx_mor na01 |
is02 : cx_mor na02 |
is12 : cx_mor na12 |
is21 : cx_mor na21 .

Ltac is_cx_mor :=
match goal with
| |- cx_mor na00 => ap is00
| |- cx_mor na11 => ap is11
| |- cx_mor na22 => ap is22
| |- cx_mor na01 => ap is01
| |- cx_mor na02 => ap is02
| |- cx_mor na12 => ap is12
| |- cx_mor na21 => ap is21
| _ => fail end.
 

Lemma cx_mor_mor : forall u, cx_mor u ->
mor (coarse_cat (R 3)) u.

Definition cx := subcategory (coarse_cat (R 3))
cx_ob cx_mor.

Ltac R_nat_discriminate :=
match goal with
| id1 : R ?X1 = R ?X2 |- _ =>
assert (discr : X1 = X2); [app R_inj | discriminate discr]
| _ => fail
end.

Lemma cx_property : subcategory_property (coarse_cat (R 3))
cx_ob cx_mor.

Lemma cx_axioms : Category.axioms cx.

Lemma ob_cx : forall x, ob cx x = cx_ob x.

Lemma mor_cx : forall u, mor cx u = cx_mor u.

Inductive in_cx_sys : E -> Prop :=
sys00 : in_cx_sys na00 |
sys11 : in_cx_sys na11 |
sys22 : in_cx_sys na22 |
sys01 : in_cx_sys na01 |
sys02 : in_cx_sys na02.

Ltac is_cx_sys :=
match goal with
| |- in_cx_sys na00 => ap sys00
| |- in_cx_sys na11 => ap sys11
| |- in_cx_sys na22 => ap sys22
| |- in_cx_sys na01 => ap sys01
| |- in_cx_sys na02 => ap sys02
| _ => fail end.

Definition cx_sys := Z (morphisms cx) in_cx_sys.

Lemma inc_cx_sys : forall u,
inc u cx_sys = in_cx_sys u.

Lemma source_null_arrow_eq : forall a b c d,
null_arrow a b = null_arrow c d -> a = c.

Lemma target_null_arrow_eq : forall a b c d,
null_arrow a b = null_arrow c d -> b = d.

Lemma cx_sys_rw : forall u,
inc u cx_sys = (in_cx_sys u &
mor cx u & cx_mor u & ~u=na12 & ~u=na21).

Lemma comp_cx : forall u v,
mor cx u -> mor cx v -> source u = target v ->
comp cx u v = na_comp u v.

Lemma id_cx : forall x,
ob cx x -> id cx x = null_arrow x x.

Lemma localizing_system_cx_sys:
localizing_system cx cx_sys.

Lemma multiplicative_system_cx_sys:
multiplicative_system cx cx_sys.

Lemma has_left_fractions_cx_sys :
has_left_fractions cx cx_sys.


Definition lf110 := lf_symbol na11 na01.

Definition lf120 := lf_symbol na12 na02.

Lemma lff110 : lf_forward lf110 = na11.

Lemma lfb110 : lf_backward lf110 = na01.

Lemma lff120 : lf_forward lf120 = na12.

Lemma lfb120 : lf_backward lf120 = na02.

Lemma is_lf_symbol_lf110 : is_lf_symbol cx cx_sys lf110.

Lemma is_lf_symbol_lf120 : is_lf_symbol cx cx_sys lf120.


Lemma lfv110 : lf_vertex lf110 = R 1.

Lemma lfv120 : lf_vertex lf120 = R 2.

Lemma under_lf110_same : forall u,
lf_under cx cx_sys lf110 u -> u = lf110.

Lemma under_lf120_same : forall u,
lf_under cx cx_sys lf120 u -> u = lf120.

Lemma lf110_different_lf120 : ~lf110 = lf120.

Lemma lf_beyond_lf110_lf120 : lf_beyond cx cx_sys lf110 lf120.

Lemma beyond_not_under_counterexample :
exists a, exists s, exists u, exists v,
(has_left_fractions a s &
is_lf_symbol a s u &
is_lf_symbol a s v &
lf_beyond a s u v &
~(exists w, (lf_under a s u w & lf_under a s v w))).

End Lf_Counterexample.

\end{verbatim}

\newpage

\section{The file {\tt infinite.v}}

\begin{verbatim}

Require Export gzloc.
Require Export cardinal.

\end{verbatim}

\subsection{The module {\tt Infinite}}

\begin{verbatim}

Module Infinite.
Export UpdateA.
Export Cardinal.

Lemma is_infinite_rw1 : forall x,
is_infinite x =
(exists f, Transformation.injective x x f &
~(Transformation.surjective x x f)).

Lemma is_infinite_rw : forall x,
is_infinite x =
(exists f, exists u, (inc u x & Transformation.injective x x f &
(forall v, inc v x -> ~ (f v = u)))).

Lemma infinite_sub_infinite : forall x,
(exists y, (is_infinite y & sub y x)) -> is_infinite x.

Lemma is_finite_natural : forall x, inc x nat ->
is_finite x.

Lemma is_infinite_nat : is_infinite nat.

Lemma is_finite_cardinality : forall x,
is_finite (cardinality x) = is_finite x.

Lemma is_infinite_cardinality : forall x,
is_infinite (cardinality x) = is_infinite x.

Definition natE : E.

Lemma is_ordinal_nat : is_ordinal nat.

Lemma inc_cardinality_nat_or_sub_nat_cardinality : forall x,
(inc (cardinality x) nat) \/ (sub nat (cardinality x)).

Lemma is_finite_inc_cardinality_nat : forall x,
is_finite x = inc (cardinality x) nat.

Lemma is_infinite_sub_nat_cardinality : forall x,
is_infinite x = sub nat (cardinality x).

Definition plus_one x := tack_on x x.

Lemma is_ordinal_plus_one : forall x,
is_ordinal x -> is_ordinal (plus_one x).

Lemma plus_one_R : forall (i:nat),
plus_one (R i) = R (i+1).

Lemma inc_plus_one_nat : forall x,
inc x nat -> inc (plus_one x) nat.

Lemma function_V_tack_on_old : forall f x y z,
Function.axioms f -> ~inc x (domain f) ->
inc z (domain f) ->
V z (tack_on f (pair x y)) = V z f.

Lemma function_V_tack_on_new : forall f x y z,
Function.axioms f -> ~inc x (domain f) ->
z = x ->
V z (tack_on f (pair x y)) = y.

Lemma are_isomorphic_tack_on : forall x y u v,
~inc y x -> ~inc v u ->
are_isomorphic x u ->
are_isomorphic (tack_on x y) (tack_on u v).

Lemma nat_are_isomorphic_eq : forall x y,
inc x nat -> inc y nat -> are_isomorphic x y ->
x = y.

Lemma nat_cardinality_refl : forall x,
inc x nat -> cardinality x = x.

Lemma cardinality_tack_on : forall x y,
is_finite x -> ~inc y x ->
cardinality (tack_on x y) = plus_one (cardinality x).

Definition number x := Bnat (cardinality x).

Lemma R_number : forall x,
is_finite x -> R (number x) = cardinality x.

Lemma cardinality_emptyset :
cardinality emptyset = emptyset.

Lemma number_emptyset : forall x,
x = emptyset -> number x = 0.

Lemma number_tack_on : forall x y,
is_finite x -> ~inc y x ->
number (tack_on x y) = (number x) + 1.

Definition take_out x y :=
complement x (singleton y).

Lemma inc_take_out : forall x y z,
inc z (take_out x y) =
(inc z x & ~z=y).

Lemma tack_on_take_out : forall x y,
inc y x ->
tack_on (take_out x y) y = x.

Lemma not_inc_take_out : forall x y,
~inc y (take_out x y).

Lemma subset_finite : forall x,
(exists y, (sub x y & is_finite y)) -> is_finite x.

Lemma is_finite_take_out : forall x y,
is_finite x -> is_finite (take_out x y).

Lemma number_take_out : forall x y,
is_finite x -> inc y x -> number (take_out x y) = (number x) -1.

Lemma number_zero_emptyset : forall x,
is_finite x -> number x = 0 -> x = emptyset.

Lemma number_gt_zero_nonempty : forall x,
is_finite x -> number x > 0 ->
(exists y, inc y x).

Lemma finite_induction: forall (p:EP) x,
(p emptyset) ->
(forall y z, p y -> p (tack_on y z)) ->
is_finite x -> p x.

Lemma image_create_tack_on : forall x y f,
Image.create (tack_on x y) f =
tack_on (Image.create x f) (f y).

Lemma image_create_emptyset : forall f,
Image.create emptyset f = emptyset.

Lemma image_finite : forall x f,
is_finite x -> is_finite (Image.create x f).


Definition is_ordinal_pair x :=
is_pair x & is_ordinal (pr1 x) & is_ordinal (pr2 x).

Definition op_lex_lt x y :=
is_ordinal_pair x &
is_ordinal_pair y &
(ordinal_lt (pr1 x) (pr1 y) \/
(pr1 x = pr1 y & ordinal_lt (pr2 x) (pr2 y))).

Definition op_lex_leq x y :=
is_ordinal_pair x &
is_ordinal_pair y &
(ordinal_lt (pr1 x) (pr1 y) \/
(pr1 x = pr1 y & ordinal_leq (pr2 x) (pr2 y))).

Lemma ordinal_leq_refl : forall x,
is_ordinal x -> ordinal_leq x x.

Lemma ordinal_lt_leq : forall x y,
ordinal_lt x y -> ordinal_leq x y.

Lemma ordinal_leq_rw : forall x y,
ordinal_leq x y =
(ordinal_lt x y \/ (x = y & is_ordinal x)).

Lemma op_lex_lt_leq : forall x y,
op_lex_lt x y -> op_lex_leq x y.

Lemma op_lex_leq_refl : forall x,
is_ordinal_pair x -> op_lex_leq x x.

Lemma op_lex_leq_rw : forall x y,
op_lex_leq x y =
(op_lex_lt x y \/ (x = y & is_ordinal_pair x)).

Lemma not_ordinal_lt_refl : forall x,
~(ordinal_lt x x).

Lemma not_op_lex_lt_refl : forall x,
~ (op_lex_lt x x).

Lemma ordinal_lt_trans : forall x z,
(exists y, (ordinal_lt x y & ordinal_lt y z)) ->
ordinal_lt x z.

Lemma op_lex_lt_trans : forall x z,
(exists y, (op_lex_lt x y & op_lex_lt y z)) ->
op_lex_lt x z.

Lemma op_lex_leq_trans : forall x z,
(exists y, (op_lex_leq x y & op_lex_leq y z)) ->
op_lex_leq x z.

Lemma not_ordinal_lt_lt : forall x y,
ordinal_lt x y -> ordinal_lt y x -> False.

Lemma not_op_lex_lt_lt : forall x y,
op_lex_lt x y -> op_lex_lt y x -> False.

Lemma op_lex_leq_leq_eq : forall x y,
op_lex_leq x y -> op_lex_leq y x -> x = y.

Lemma op_lex_dichotomy : forall x y,
is_ordinal_pair x -> is_ordinal_pair y ->
(op_lex_lt x y \/ op_lex_lt y x \/ x = y).

Definition least_ordinal z := sow (fun x => inc x z).

Lemma is_ordinal_least_ordinal: forall z,
is_ordinal (least_ordinal z).

Lemma inc_least_ordinal_ex : forall z,
(exists y, (inc y z & is_ordinal y)) ->
inc (least_ordinal z) z.

Definition set_of_ordinals z :=
(forall y, inc y z -> is_ordinal y).

Lemma inc_least_ordinal : forall z,
set_of_ordinals z -> nonempty z ->
inc (least_ordinal z) z.

Lemma ordinal_leq_least_ordinal : forall y z,
inc y z -> is_ordinal y ->
ordinal_leq (least_ordinal z) y.

Definition set_of_ordinal_pairs z :=
(forall y, inc y z -> is_ordinal_pair y).

Lemma set_of_ordinal_pairs_relation : forall z,
set_of_ordinal_pairs z -> is_relation z.

Lemma set_of_ordinals_domain : forall z,
set_of_ordinal_pairs z ->
set_of_ordinals (domain z).

Definition column z x:=
range (Z z (fun y => pr1 y = x)).

Lemma inc_column : forall z x y,
is_relation z ->
inc y (column z x) = inc (pair x y) z.

Lemma set_of_ordinals_column : forall z x,
set_of_ordinal_pairs z -> set_of_ordinals (column z x).

Lemma nonempty_domain : forall z,
nonempty z -> nonempty (domain z).

Lemma nonempty_column : forall z x,
is_relation z ->
inc x (domain z) -> nonempty (column z x).

Definition op_lex_least z :=
pair (least_ordinal (domain z))
(least_ordinal (column z (least_ordinal (domain z)))).

Lemma is_ordinal_pair_op_lex_least : forall z,
is_ordinal_pair (op_lex_least z).

Lemma inc_op_lex_least : forall z,
set_of_ordinal_pairs z -> nonempty z ->
inc (op_lex_least z) z.

Lemma pr1_op_lex_least : forall z,
pr1 (op_lex_least z) = least_ordinal (domain z).

Lemma pr2_op_lex_least : forall z,
pr2 (op_lex_least z) =
least_ordinal (column z (least_ordinal (domain z))).

Lemma op_lex_leq_op_lex_least : forall z y,
set_of_ordinal_pairs z -> inc y z ->
op_lex_leq (op_lex_least z) y.


Definition op_lex z := Order_notation.create z op_lex_leq.

Lemma U_op_lex : forall z, U (op_lex z) = z.

Lemma leq_op_lex : forall z u v:E,
leq (op_lex z) u v = (inc u z & inc v z & op_lex_leq u v).

Lemma op_lex_axioms : forall z:E,
set_of_ordinal_pairs z ->
Order.Definitions.axioms (op_lex z).

Lemma is_linear_op_lex : forall z:E,
set_of_ordinal_pairs z ->
is_linear (op_lex z).

Lemma is_well_ordered_op_lex : forall z,
set_of_ordinal_pairs z ->
is_well_ordered (op_lex z).

Definition prod_reflexive z :=
are_isomorphic z (Cartesian.product z z).

Lemma prod_reflexive_invariant : forall y z,
are_isomorphic y z -> prod_reflexive y -> prod_reflexive z.

Definition sub_prod_reflexive z :=
forall y, sub y z -> is_infinite y -> prod_reflexive y.

Lemma sub_prod_reflexive_rw : forall z,
sub_prod_reflexive z =
(forall y, iso_sub y z -> is_infinite y -> prod_reflexive y).

Lemma sub_prod_reflexive_invariant : forall y z,
iso_sub y z -> sub_prod_reflexive z -> sub_prod_reflexive y.

Lemma sub_prod_reflexive_prod_reflexive : forall z,
is_infinite z -> sub_prod_reflexive z -> prod_reflexive z.

Lemma finite_sub_prod_reflexive : forall z,
is_finite z -> sub_prod_reflexive z.

Definition op_triangle z :=
Z (Cartesian.product z z) (fun y => (ordinal_leq (pr2 y) (pr1 y))).

Lemma inc_op_triangle : forall y z,
is_ordinal z -> inc y (op_triangle z) =
(is_ordinal_pair y & ordinal_lt (pr1 y) z &
ordinal_leq (pr2 y) (pr1 y)).

Lemma set_of_ordinal_pairs_op_triangle : forall z,
is_ordinal z -> set_of_ordinal_pairs (op_triangle z).

Lemma iso_sub_op_triangle : forall z,
is_ordinal z -> iso_sub z (op_triangle z).

Lemma sub_op_triangle_product : forall z,
sub (op_triangle z) (product z z).

Lemma same_cardinality_op_triangle : forall z,
is_infinite z -> is_ordinal z -> sub_prod_reflexive z ->
cardinality (op_triangle z) = cardinality z.

Lemma op_triangle_sub : forall y z,
sub y z -> sub (op_triangle y) (op_triangle z).

Definition op_triangle_lex z := op_lex (op_triangle z).

Lemma U_op_triangle_lex : forall z,
U (op_triangle_lex z) = op_triangle z.

Lemma leq_op_triangle_lex : forall z u v,
leq (op_triangle_lex z) u v = (inc u (op_triangle z) &
inc v (op_triangle z) & op_lex_leq u v).

Lemma is_well_ordered_op_triangle_lex : forall z,
is_ordinal z -> is_well_ordered (op_triangle_lex z).

Definition is_cardinal z :=
cardinality z = z.

Lemma is_cardinal_rw : forall z,
is_cardinal z = (is_ordinal z & (forall y, ordinal_lt y z ->
ordinal_lt (cardinality y) (cardinality z))).

Lemma is_cardinal_rw2 : forall z,
is_cardinal z = (is_ordinal z & (forall y, is_ordinal y ->
iso_sub z y -> ordinal_leq z y)).

Lemma is_cardinal_cardinality : forall z,
is_cardinal (cardinality z).

Lemma infinite_cardinality_tack_on : forall x y,
is_infinite x -> cardinality (tack_on x y) = cardinality x.

Lemma cardinality_lt_ordinal_lt : forall y z,
is_ordinal y -> is_ordinal z ->
ordinal_lt (cardinality y) (cardinality z) -> ordinal_lt y z.

Lemma ord_card_leq_criterion : forall y z,
is_ordinal y -> is_cardinal z ->
(forall x, inc x y -> ordinal_lt x z) ->
ordinal_leq y z.

Lemma inc_punctured_downward_subset : forall a x y,
Order.Definitions.axioms a -> inc x (U a) ->
inc y (punctured_downward_subset a x) = lt a y x.

Lemma punctured_downward_subset_wo_avatar : forall a x,
is_well_ordered a -> inc x (U a) ->
Transformation.bijective (punctured_downward_subset a x)
(wo_avatar a x) (wo_avatar a).

Lemma cardinality_punctured_downward_subset : forall a x,
is_well_ordered a -> inc x (U a) ->
(cardinality (punctured_downward_subset a x) =
cardinality (wo_avatar a x)).
 

Lemma card_leq_criterion : forall a z,
is_well_ordered a -> is_cardinal z ->
(forall x, inc x (U a) ->
ordinal_lt (cardinality (punctured_downward_subset a x)) z) ->
ordinal_leq (cardinality (U a)) z.

Lemma sub_punctured_downward_cone_op_triangle : forall y z,
is_ordinal z -> inc y (op_triangle z) ->
sub (punctured_downward_subset (op_triangle_lex z) y)
(op_triangle (tack_on (pr1 y) (pr1 y))).

Definition sub_prod_reflexive_below z :=
is_cardinal z & is_infinite z &
(forall y, ordinal_lt y z -> sub_prod_reflexive y).

Lemma ordinal_leq_lt_trans : forall x y,
(exists z, (ordinal_leq x z & ordinal_lt z y)) ->
ordinal_lt x y.

Lemma is_finite_union2 : forall a b, is_finite a ->
is_finite b -> is_finite (union2 a b).

Lemma product_singleton_isomorphic : forall a x,
are_isomorphic a (product a (singleton x)).

Lemma is_finite_product : forall a b,
is_finite a -> is_finite b -> is_finite (product a b).

Lemma infinite_tack_on_isomorphic: forall z x,
is_infinite z -> are_isomorphic z (tack_on z x).

Lemma sprb_op_triangle_leq : forall z,
sub_prod_reflexive_below z ->
ordinal_leq (cardinality (op_triangle z)) z.

Lemma inc_tack_on_infinite_cardinal : forall z x,
is_cardinal z -> is_infinite z -> inc x z -> inc (tack_on x x) z.

Lemma tack_on_injective : forall a b,
is_ordinal a -> is_ordinal b ->
tack_on a a = tack_on b b -> a = b.

Lemma iso_sub_double_op_triangle : forall z,
is_cardinal z -> is_infinite z ->
iso_sub (product z (R 2)) (op_triangle z).

Lemma finite_ordinal_cardinal : forall x,
is_ordinal x -> is_finite x -> is_cardinal x.

Lemma ordinal_lt_finite_infinite : forall x y,
is_ordinal x -> is_ordinal y -> is_finite x
-> is_infinite y -> ordinal_lt x y.

Lemma dumb_sub_product : forall z,
is_ordinal z -> is_infinite z ->
sub (product z (R 2)) (product z z).

Definition lim_sprb z :=
sub_prod_reflexive_below (cardinality z) \/
sub_prod_reflexive z .

Lemma cardinality_double : forall z,
is_ordinal z -> is_infinite z ->
lim_sprb z ->
cardinality (product z (R 2)) = cardinality z.

Lemma cardinal_leq_image : forall z f,
ordinal_leq (cardinality (Image.create z f)) (cardinality z).

Lemma cardinal_leq_iso_sub : forall a b,
ordinal_leq (cardinality a) (cardinality b) ->
iso_sub a b.

Lemma ordinal_leq_cardinality_union2_double : forall a b z,
ordinal_leq (cardinality a) (cardinality z) ->
ordinal_leq (cardinality b) (cardinality z) ->
ordinal_leq (cardinality (union2 a b))
(cardinality (product z (R 2))).

Lemma lim_sprb_cardinal_union2 : forall a b z,
lim_sprb z -> is_infinite z ->
ordinal_leq (cardinality a) (cardinality z) ->
ordinal_leq (cardinality b) (cardinality z) ->
ordinal_leq (cardinality (union2 a b)) (cardinality z).

Definition reflect_pair x := pair (pr2 x) (pr1 x).

Definition reflected_op_triangle z :=
Image.create (op_triangle z) reflect_pair.

Lemma cardinal_leq_reflected_op_triangle : forall z,
ordinal_leq (cardinality (reflected_op_triangle z))
(cardinality (op_triangle z)).

Lemma inc_ordinal_prod_itself_or : forall z x,
is_ordinal z ->
inc x (product z z) =
(inc x (op_triangle z) \/ inc x (reflected_op_triangle z)).

Lemma ordinal_prod_itself_decomposition : forall z,
is_ordinal z ->
product z z = union2 (op_triangle z) (reflected_op_triangle z).

Lemma lim_sprb_cardinal_product : forall z,
lim_sprb z -> is_infinite z ->
ordinal_leq (cardinality (product z z)) (cardinality z).

Lemma lim_sprb_sub_prod_reflexive : forall z,
lim_sprb z -> sub_prod_reflexive z.

Lemma ordinal_lt_leq_trans : forall x z,
(exists y, ordinal_lt x y & ordinal_leq y z) ->
ordinal_lt x z.

Lemma all_ordinals_sub_prod_reflexive : forall z,
is_ordinal z -> sub_prod_reflexive z.

Lemma all_sub_prod_reflexive : forall z,
sub_prod_reflexive z.

Lemma are_isomorphic_product : forall z,
is_infinite z -> are_isomorphic z (product z z).

Lemma cardinality_product_infinite : forall z,
is_infinite z -> cardinality (product z z) = cardinality z.

Lemma cardinality_union2_infinite : forall a b z,
is_infinite z ->
ordinal_leq (cardinality a) (cardinality z) ->
ordinal_leq (cardinality b) (cardinality z) ->
ordinal_leq (cardinality (union2 a b)) (cardinality z).

Lemma iso_sub_union2 : forall a b z,
is_infinite z -> iso_sub a z ->iso_sub b z ->
iso_sub (union2 a b) z.


Lemma russell : forall z, ~(iso_sub (powerset z) z).

Lemma cardinal_lt_powerset : forall z,
ordinal_lt (cardinality z) (cardinality (powerset z)).

End Infinite.

\end{verbatim}

\end{document}